\documentclass{article}
\pdfoutput=1
\usepackage{mathptmx}
\usepackage{amsmath, amsbsy, amsfonts}
\usepackage{epsfig}
\usepackage{color}
\usepackage{url}
\usepackage{authblk}
\newcommand{\gr}{\mbox{$\nabla$}}

\newcommand{\dv}{\mbox{$\mbox{div}$}}

\newcommand{\Ib}{\mathbb{I}}
\newcommand{\Kb}{\mathbb{K}}
\newcommand{\Ab}{\mathbb{A}}
\newcommand{\xb}{\mbox{\bf x}}

\newcommand{\der}[2]{\frac{\partial #1}{\partial #2}}
\newcommand{\dert}[1]{\frac{\partial #1}{\partial t}}

\newcommand{\ql}{\mbox{$\mathbf{q}_l$}}
\newcommand{\qg}{\mbox{$\mathbf{q}_g$}}
\newcommand{\gb}{\mbox{$\mathbf{g}$}}
\newcommand{\jlw}{\mbox{$\mathbf{j}_l^w$}}
\newcommand{\jlh}{\mbox{$\mathbf{j}_l^h$}}

\newcommand{\Htwo}{$\text{H}_2$}
\newtheorem{remark}{Remark}
\newtheorem{acknowledgements}{Acknowledgements}
\def\keywords{\vspace{.5em}
{\textit{Keywords}:\,\relax%
}}

\def\Ind{\ensuremath{\mbox{1\hspace{-.25em}l}}}
\author[1]{Alain Bourgeat}  
\author[2]{Mladen Jurak}
\author[3]{Farid Sma\"{\i}}
\affil[1]{Universit\'{e} de Lyon, Universit\'{e}
Lyon1, CNRS UMR 5208 Institut Camille Jordan, F - 69200 Villeurbanne
Cedex, France}
\affil[2]{Department of mathematics,
University of Zagreb,Bijenicka 30, Zagreb, Croatia,
}
\affil[3]{Universit\'{e} de Lyon, Universit\'{e} Lyon1,
CNRS UMR 5208 Institut Camille Jordan, F - 69200 Villeurbanne Cedex,
France } 
\date{\today}
\title{Modelling and Numerical Simulation of gas migration in a nuclear waste repository}
\begin{document}

\maketitle
\begin{abstract}
We present a compositional compressible two-phase, liquid and gas,
flow model for numerical simulations of hydrogen migration in deep
geological radioactive waste repository. This model includes
capillary effects and the gas diffusivity.  The choice of the main
variables in this model, Total or Dissolved Hydrogen Mass
Concentration and Liquid Pressure, leads to a unique and consistent
formulation of the  gas phase appearance and disappearance. After
introducing this model, we show computational evidences of its
adequacy to simulate  gas phase appearance and disappearance in
different situations typical of underground radioactive waste
repository.
\end{abstract}
\keywords{Two-phase flow, compositional flow,  porous medium, underground
nuclear waste management}
 \tableofcontents

\section{Introduction}

The simultaneous flow of immiscible fluids in porous media occurs in
a wide variety of applications. The most concentrated research in
the field of multiphase flows over the past four decades has focused
on  unsaturated groundwater flows, and flows in underground
petroleum reservoirs. Most recently, multiphase flows have generated
serious interest among engineers concerned with  deep geological
repository for radioactive waste. There is growing awareness that
the effect of hydrogen gas generation, due to anaerobic corrosion of
the steel engineered barriers of radioactive waste packages(carbon
steel overpacks and stainless steel envelopes), can affect all the
functions allocated to the canisters or to the buffers and the
backfill.
 The host rock safety function may even be
threaten by overpressurisation leading to opening fractures in the
host rock and inducing groundwater flow and transport of
radionuclides outside of the waste site boundaries.

 Equations governing this type of flow in porous media are inherently nonlinear,
and the geometries and material properties characterizing many
situations in many applications ( petroleum reservoir, gas storage,
waste repository),  can be quite irregular and contrasted. As a
result of all these difficulties, numerical simulation often offers
the only viable approach to modelling multiphase porous-media flows . In
nuclear waste management, the migration of gas through the near
field environment and the host rock,
 involves two components, water and pure hydrogen \Htwo; and
two phases "liquid" and "gas".  Our ability to understand and
predict underground gas migration  is crucial to the design and to
assessing the performance of reliable nuclear waste storages. This
is a fairly new frontier in multiphase porous-media flows, and again
the inherent complexity of the physics leads to governing equations
for which the only practical way to produce solutions may be
numerical simulation.

 This paper addresses one of the outstanding
physical and
 mathematical problems in
multiphase flow simulation:  the appearance and disappearance of one
of the phases, leading to the degeneracy of the equations satisfied
by the saturation. In order to overcome this difficulty, we  discuss
a formulation
 based on  variables  which doesn't degenerate and hence could be used as
 an unique formulation for both situations, liquid saturated and unsaturated.
We will demonstrate through four numerical tests, the ability of
this new formulation to actually cope with the appearance or/and
disappearance  of one phase in simple, typical but challenging
situations, like the ones we met in underground radioactive waste
repository simulations.

\section{ Modeling Physical Assumptions}
We consider herein a
 porous medium saturated with a fluid composed of 2 phases,
 \textit{liquid} and \textit{gas}, and two
 components. According to the application we have in mind, we
 consider the
 fluid as a mixture of two
 components: water (only liquid) and hydrogen (\Htwo, mostly gas) or any gas with similar thermodynamical
 properties. In the following, for sake of simplicity we will call
 \textit{hydrogen}  the non-water component and use indices $w$ and $h$ for the  \textit{water} and the
 \textit{hydrogen} components.

 We neglect the  water vaporization since, in underground formations with high water pressure, the water vapor does not contribute
 significantly to the gas phase pressure .
 The water component is incompressible while the gas phase follow
  the ideal gas law.  The whole fluid system is in thermal equilibrium and
  the porous medium is rigid, meaning that  the
 porosity $\Phi$ is only a function of the space variable
 $\Phi=\Phi(\xb)$;
moreover, since hydrogen is  highly  diffusive  we  include the
dissolved
 hydrogen diffusion
 in the liquid phase .

The two phases are denoted by indices, $l$ for liquid, and $g$ for
gas. Associated to each phase $\alpha\in\{l,g\}$, we have, in the
porous medium,  the phase pressures $p_\alpha$, the phase
saturations $S_\alpha$, the phase mass densities $\rho_\alpha$ and
the phase volumetric flow rates $\mathbf{q}_\alpha$.
The phase volumetric flow rates are given by the {\em  Darcy-Muskat
law}:
\begin{equation}
 \ql = - \Kb(\xb)\lambda_l(S_l)
          \left( \gr p_l -\rho_l \gb\right),\quad
          \qg  = - \Kb(\xb)\lambda_g(S_g)
          \left( \gr p_g -\rho_g \gb\right),\label{darcy}
\end{equation}
where $ \Kb(\xb)$ is the absolute permeability tensor,
$\lambda_\alpha(S_\alpha)$ is  the $\alpha-$phase relative mobility
function, and $\gb$ is the gravity acceleration;
 $S_\alpha$ is the reduced $\alpha-$phase saturation and then
 satisfies:
\begin{equation}
S_l + S_g = 1 .    \label{saturations}
\end{equation}
Pressures are connected through a given  {\em capillary pressure
law}:
\begin{equation}
 p_c(S_g) = p_g -p_l.
 \label{capillary}
\end{equation}
From definition (\ref{capillary}) we notice that $p_c$ is a strictly
increasing function of gas saturation, $ p_c'(S_g) > 0 $, leading to
a {\em capillary constraint}:
\begin{equation}
 p_g  > p_l + p_c(0),
 \label{capillaryconstraint}
\end{equation}
 where $p_c(0)\geq 0$  is the capillary curve entry pressure ( see Figure \ref{P-C}).

Since the liquid phase could be composed of
water and dissolved hydrogen, we need to introduce the water mass
concentration $\rho_l^w$ in the liquid phase, and the hydrogen mass
concentration $\rho_l^h$ in the liquid phase. Note that the upper
index is the component index, and the lower one denotes the phase.
We have, then
\begin{equation}
\rho_l = \rho_l^w + \rho_l^h.
\label{partial-densities}
\end{equation}
As said before, in the gas phase,  we neglect the water vaporization
and we use the ideal gas law:
\begin{equation}
\rho_g = C_v p_g, \label{ideal}
\end{equation}
with $C_v = M^h/(RT)$, where $T$ is the temperature, $R$ the
universal gas constant and  $M^h$ the hydrogen molar mass. Mass
conservation for each component leads to the following differential
equations:
\begin{align}
 \Phi\dert{}\left(  S_l \rho_l^w  \right)   + \dv \left(  \rho_l^w  \ql   +\jlw \right) = {\cal F}^w, \label{eq-1} \\
\Phi\dert{}\left(  S_l \rho_l^h  + S_g  \rho_g\right) + \dv \left(  \rho_l^h \ql +  \rho_g  \qg +\jlh \right)
= {\cal F}^h,\label{eq-2}
\end{align}
where the phase flow velocities, $\ql$ and $\qg$, are given by the
Darcy-Muskat law (\ref{darcy}), ${\cal F}^c$ are the source terms,
and $\mathbf{j}_l^c$, $c\in\{w,h\}$ are the $c-$component diffusive
flux in liquid  phase, as defined later in (\ref{diff-fluxes}).

Assuming  water incompressibility  and that the liquid  volume is
independent of the dissolved hydrogen concentration, we may assume
the water component concentration in the liquid phase to be
constant, i.e.:
\begin{equation}
    \rho_l^w = \rho_w^{std},
\label{const-w}
\end{equation}
where $\rho_w^{std}$ is the standard water density.

 The assumption of hydrogen thermodynamical equilibrium in both phases leads to equal chemical potentials in each
 phase:
 $\mu_{g}^{h}(T,p_{g},X_{l}^{h})=\mu_{l}^{h}(T,p_{l},X_{l}^{h})$.
Assuming  that in the gas phase there is only the hydrogen component
and no water, leads to $X_{g}^{h}=1$; and then, from the above
chemical potentials equality,  we have a relationship
$p_{g}=F(T,p_{l},X_{l}^{h})$. Assuming that the liquid pressure
influence could be neglected  in the pressure range considered
herein and using the hydrogen low solubility, $\rho_l^h \ll
\rho_l^w=\rho_w^{std}$, we may then linearize the relationship
between $p_g$ and $X_l^h$, and finally obtain the {\em Henry's law}
$p_g = K^h X_l^h$, where $K^h$  is specific to the mixture
water/hydrogen and depends only on the temperature $T$.
 Furthermore, using (\ref{const-w}) and the hydrogen low
solubility,
 the molar  fraction, $X_l^h$, reduces to $\frac{\rho_l^{h} M^w }{\rho_w^{std} M^h}$
 (see eqs.(9)-(11) in \cite{TA02-BJS}) ~and  the {\em Henry law} can be written as
\begin{equation}
 \rho_l^h= C_h p_g,   \label{Henry}
\end{equation}
 where $C_h = H M^h = \rho_w^{std} M^h/(M^w
K^h)$; $H$ is called the Henry law constant and is also depending
only on the temperature.
\begin{figure}
\centering
\includegraphics[angle=0, height=1.6in]{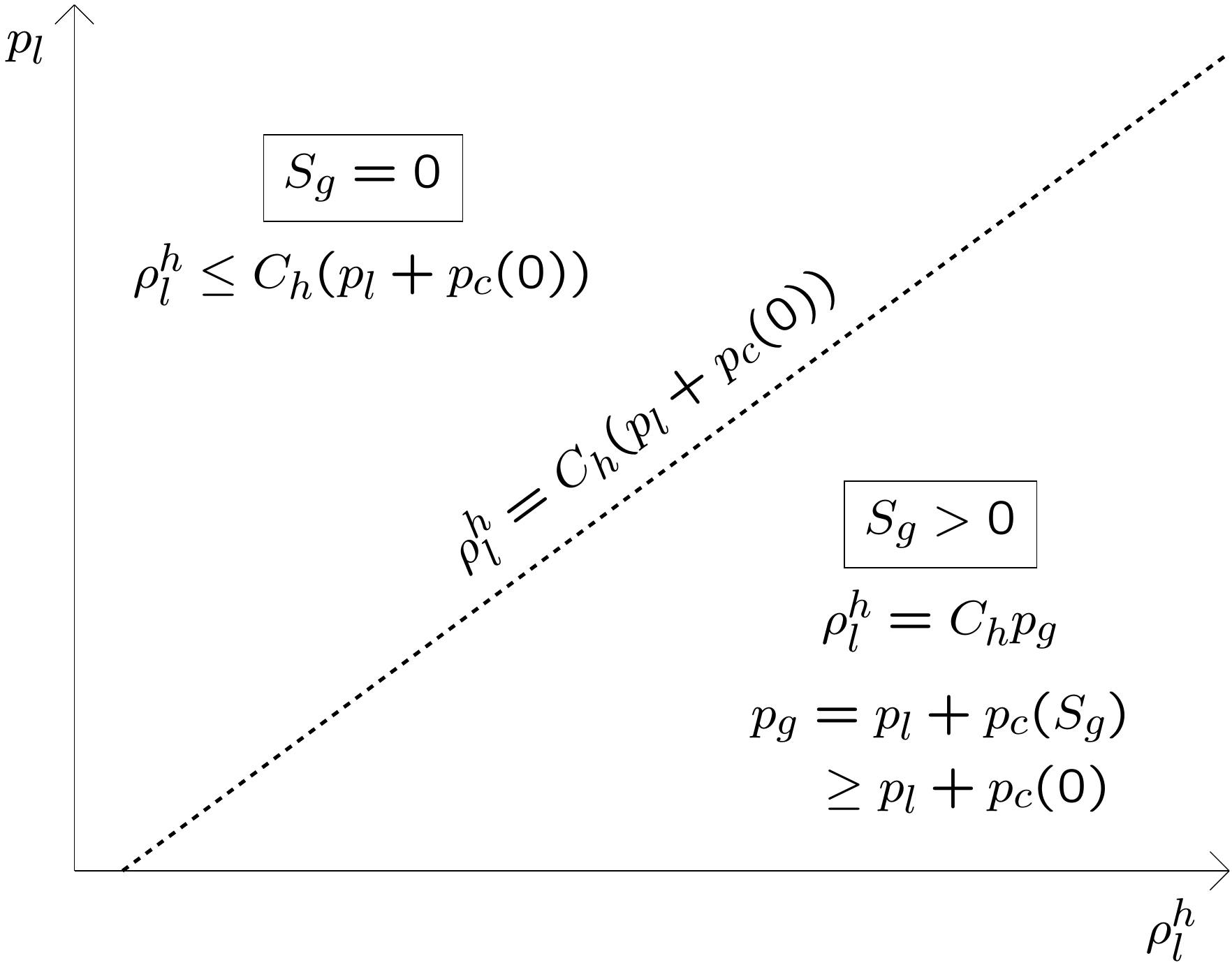}
\caption{Phase diagram:  Henry's law; localization of the liquid
saturated $S_g=0$ and  unsaturated $S_g>0$  states } \label{P-D}
\end{figure}
        \begin{remark}
        \label{Rem:1}
On the one hand the gas pressure obey the Capillary pressure law
(\ref{capillary}) with the constraint (\ref{capillaryconstraint}),
but on the other hand it should also satisfy the local
thermodynamical equilibrium and obey the Henry law (\ref{Henry}).
More precisely if there are two phases, i.e. if the concentration,
$\rho_l^h$, is sufficiently high to have a gas phase
appearance($S_g>0$) , we have from (\ref{Henry}) and
(\ref{capillary}) :
    \begin{equation}
\rho_l^h = C_h(p_l +p_c(S_g)). \label{threshold}
    \end{equation}
 Moreover, $S_g>0$ with the constraint (\ref{capillaryconstraint})
 and the Henry's law
(\ref{Henry}), gives the constraint:
    \begin{equation}
{\rho_l^h}> C_h(p_l +p_c(0)). \label{seuil}
\end{equation}
But if the concentration, $\rho_l^h$, is smaller than a certain
concentration threshold (see Figure \ref{P-D}), then there is only
the liquid phase (no gas phase, $S_g=0$), and none of all the
relationships (\ref{capillary}) or (\ref{seuil}), connected to
capillary equilibrium, applies anymore; we have only $S_g=0$, with
$\rho_l^h\leq C_h p_g$.\\
There is then a concentration threshold line, corresponding to
$\rho_l^h = C_h(p_l +p_c(0))$ in the phase diagram (Fig.\ref{P-D}),
separating the one phase (liquid saturated) region from the two
phase (liquid unsaturated) region.
        \label{rem_nonequ}
        \end{remark}

The existence of a concentration threshold line can also be written
as an unilateral condition:
\begin{align*}
    0\leq S_g \leq 1,\quad 0\leq \rho_l^h \leq C_h p_g,\quad S_g (C_h p_g -\rho_l^h) =
    0;
\end{align*}
which could be then used (see \cite{Jaffre}) for designing  a
numerical scheme based on approximating a variational equation.
 \hfill $\square$ \\

The diffusive fluxes in the liquid phase are given by the Fick law
applied to $X_l^w$ and to $X_l^h$, the water component and the
hydrogen component molar fractions (see eqs.(12) and (13) in
\cite{TA02-BJS}).
 Using the same kind of approximation as in the
Henry law, based on the hydrogen low solubility,
  we obtain, for the diffusive fluxes in this binary mixture (see Remark 2 and Remark 3
 in \cite{TA02-BJS}):
\begin{align}
\mathbf{j}_l^h =  -\Phi S_l D  \nabla \rho_l^h,\quad \mathbf{j}_l^w = -\mathbf{j}_l^h,
\label{diff-fluxes}
\end{align}
where $D
$ is the hydrogen  molecular diffusion coefficient in the
liquid phase, corrected by the tortuosity of the porous medium.

If both liquid and gas phases exist, ($S_g \neq 0$),~the porous
media is said \textit{liquid unsaturated} and  the transport model
for the liquid-gas system can be now written as:
\begin{align}
  \Phi\rho_{w}^{std}\dert{S_l} +
  \dv \left(\rho_{w}^{std} \ql  - \jlh \right)={\cal F}^w, \label{consSl}\\
  \Phi\dert{}( S_l \rho_l^h   + C_v p_g S_g  )
  +\dv \Big(\rho_l^h \ql +  C_v p_g \qg  +  \jlh\Big)
  = {\cal F}^h , \label{h2-mass-cons}\\
  \ql = -\Kb\lambda_l(S_l)
  \left( \gr p_l - (\rho_{w}^{std} + \rho_l^h) \gb\right),
  \quad
  \qg = - \Kb\lambda_g(S_g) \left( \gr p_g -C_v p_g \gb\right),
  \label{darcy-vel}\\
   \jlh = -  \Phi S_l D \gr \rho_l^h.  \label{eq:5}
\end{align}
But in the  \textit{liquid saturated} regions, where the gas phase
doesn't appear, $S_l=1$ or $S_g =0$, the system
(\ref{consSl})--(\ref{eq:5}) degenerates to:
 \begin{align}
  & \dv \left( \rho_{w}^{std}  \ql -\jlh  \right)={\cal F}^w, \quad
  \Phi\dert{\rho_l^h}  + \dv \Big(
  \rho_l^h\ql  +\jlh \Big) = {\cal F}^h ;\label{un:eq:12}\\
  \ql &= -   \Kb\lambda_l(1)
            \left( \gr p_l - (\rho_{w}^{std} +\rho_l^h ) \gb\right), \quad
  \jlh = -  \Phi  D \gr\rho_l^h .
   \label{un:eq:34}
   \end{align}

\section{Liquid Saturated/Unsaturated state; a general formulation }

A typical choice for the two primary unknowns, in modeling
immiscible two-phase flow,  is the saturation and one of the phases
pressure, for example $S_g$ and $p_l$. But as seen above, in
(\ref{consSl})--(\ref{un:eq:34}), this set of unknowns obviously
cannot describe the flow in a liquid saturated region, where there
is only one phase, and cannot take in account the gas dissolution
since then the dissolved gas concentration, $\rho_l^h$, becomes an
independent unknown.

\subsection{Modeling based on Total hydrogen concentration, $\rho_{tot}^{h}$ }

To solve this problem,  instead of using the gas saturation $S_g$
we have
 proposed, in \cite{TA02-BJS}, to use $\rho_{tot}^{h}$, the
total hydrogen mass concentration, defined as:
\begin{align}
\rho_{tot}^{h}= S_l \rho_l^h+S_g\rho_g^h.
    \label{def-RHO}
\end{align}
Defining
\begin{align}
  a(S_g) = C_h (1-S_g)   + C_v S_g  \in [ C_h, C_v];
  \label{a:def}
\end{align}
with
\begin{align}
  \quad a'(S_g) = C_v - C_h =C_{\Delta} > 0,
  \label{a:pos}
\end{align}
since $C_v > C_h$, from the assumption of weak solubility; we may
then rewrite the total hydrogen mass concentration,~$
\rho_{tot}^{h}$, defined in (\ref{def-RHO}), as:

\begin{align}
    \rho_{tot}^{h} = \begin{cases}
        a(S_g) (p_l + p_c(S_g) ) & \text{if } S_g >0\\
        \rho_l^h & \text{if } S _g = 0.
    \end{cases} \label{def-X}
\end{align}

With  this new set of  unknowns, $\rho_{tot}^{h}$ and $p_l$, the two
systems of equations (\ref{consSl})--(\ref{eq:5}) and
(\ref{un:eq:12})--(\ref{un:eq:34}) now reduce to a single system of
equations:
\begin{align}
  \Phi\rho_{w}^{std}\dert{S_l} -
  \dv &\left(\rho_{w}^{std} \Kb\lambda_l(S_l)
  \left( \gr p_l - (\rho_{w}^{std} + \rho_l^h) \gb\right)  - \Phi S_l D \gr \rho_l^h \right)={\cal F}^w, \label{consSSl}\\
  \Phi\dert{\rho_{tot}^h}
  -\dv &\Big(\rho_l^h \Kb\lambda_l(S_l)
  \left( \gr p_l - (\rho_{w}^{std} + \rho_l^h) \gb\right) \nonumber\\
  &+  C_v p_g \Kb\lambda_g(S_g) \left( \gr p_l  +\gr p_c(S_g) -C_v p_g \gb\right)  +  \Phi S_l D \gr \rho_l^h\Big)
  = {\cal F}^h. \label{h2-mass-conss}
  \end{align}
 Now, if we want to study the mathematical properties of the operators in this  system of
 equations, we should develop the above system of equations
 using   $S_g= S_g(p_l,\rho_{tot}^{h})$,~
  $S_l=1-S_g=S_l(p_l,~\rho_{tot}^{h})$, and
$\rho_l^h=\rho_l^h(p_l,\rho_{tot}^{h})$,  with
\begin{align}
   \der{S_g}{p_l} =  -\frac{a(S_g)^2  \Ind_{\{\rho_{tot}^{h} > C_h (p_l +p_c(0))\}}}{C_{\Delta} \rho_{tot}^{h}
   +a(S_g)^2 p_c'(S_g)},\quad
   \der{S_g}{\rho_{tot}^{h}} =  \frac{a(S_g)  \Ind_{\{\rho_{tot}^{h} > C_h (p_l +p_c(0))\}}}{C_{\Delta} \rho_{tot}^{h}
   +a(S_g)^2 p_c'(S_g)},\label{pos:def:S}
\end{align}
where $\Ind_{\{\rho_{tot}^{h} > C_h (p_l +p_c(0))\}}$ is the
characteristic function of the set $\{ \rho_{tot}^{h} > C_h (p_l
+p_c(0)) \}$. As noted in section 2.5 in \cite{TA02-BJS}, we have
$\partial S_g/\partial p_l \leq 0$ and $\partial S_g/\partial
\rho_{tot}^{h} > 0$,  when the gas phase is present.
 Then  the system (\ref{consSl})--(\ref{h2-mass-cons}) can be written :
\begin{align}
    -\Phi\rho_w^{std}\der{S_g}{p_l}\dert{p_l} &
 - \dv \left(  \Ab^{1,1}\gr p_l + \Ab^{1,2}\gr \rho_{tot}^{h} + B_1 \Kb\gb\right)
  -\Phi\rho_w^{std}\der{S_g}{\rho_{tot}^{h}}\dert{\rho_{tot}^{h}}    ={\cal F}^w
  \label{eq:1-3}\\
  \Phi\dert{\rho_{tot}^{h}} & - \dv \Big( \Ab^{2,1}\gr p_l +\Ab^{2,2}\gr \rho_{tot}^{h} +B_2 \Kb\gb\Big)
      = {\cal F}^h.\label{eq:2-3}
\end{align}
Where the coefficients are defined by:
\begin{align}
    \Ab^{1,1}(p_l,\rho_{tot}^{h}) =& \lambda_l(S_g)\rho_{w}^{std}\Kb - {\Phi (1-S_g) } D C_h N \Ib,\label{coef:in:1}\\
     \Ab^{1,2}(p_l,\rho_{tot}^{h}) =& - {\Phi (1-S_g)  }\frac{1-N}{a(S_g)} D C_h\Ib, \\
    \Ab^{2,1}(p_l,\rho_{tot}^{h}) =& (\lambda_l(S_g)\rho_l^h +\lambda_g(S_g) C_v p_g N )\Kb +{\Phi  (1-S_g)} D C_h N\Ib,   \\
    \Ab^{2,2}(p_l,\rho_{tot}^{h}) =& \lambda_g(S_g)\frac{1-N}{a(S_g)}C_v p_g  \Kb +{\Phi(1-S_g)}\frac{1-N}{a(S_g)} D C_h\Ib, \\
     B_1(p_l,\rho_{tot}^{h}) = & -   \lambda_l(S_g)\rho_{w}^{std}  [\rho_{w}^{std} + \rho_{l}^{h}] ,\\
    B_2(p_l,\rho_{tot}^{h}) = & -(\lambda_l(S_g)  \rho_l^h [\rho_{w}^{std} +  \rho_{l}^{h}]
                 + \lambda_g(S_g)C_v^2  p_g^2) ;  \label{coef:in:2}
\end{align}
with $\Ib$ denoting the identity matrix and with the  auxiliary
functions
\begin{align}
N(p_l,\rho_{tot}^{h}) =  \frac{C_{\Delta} \rho_{tot}^{h}}{C_{\Delta}
\rho_{tot}^{h} +a(S_g)^2 p_c'(S_g)}\Ind_{\{\rho_{tot}^{h} >
C_h(p_l +p_c(0))\}}\in [0,1),\\
 \rho_l^h(p_l,\rho_{tot}^{h})  = \min(C_h  p_g(p_l,\rho_{tot}^{h}),\rho_{tot}^{h}),\quad
 p_g(p_l,\rho_{tot}^{h})  = p_l+p_c(S_g(p_l,\rho_{tot}^{h})).
\end{align}

We should notice first that equation (\ref{eq:2-3}) is uniformly
parabolic in the presence of capillarity and diffusion; but if
capillarity and diffusion are neglected, this same equation  becomes
a pure hyperbolic transport equation (see sec. 2.6 in
\cite{TA02-BJS}). Then, if we sum equations (\ref{eq:1-3}) and
(\ref{eq:2-3}) we obtain a uniformly parabolic/elliptic equation,
which is parabolic in the unsaturated (two-phases) region and
elliptic in the liquid saturated (one-phase) region.
\begin {remark}
Simulations presented in sec. 3.2 in \cite{TA02-BJS} show that this
last model with these variables, $\rho_{tot}^{h}$, the total
hydrogen mass concentration, and $p_l$, the liquid phase  pressure,
could easily handle phase transitions (appearance and disappearance
of the gas phase) in two-phase partially miscible flows.
 However, in equations (\ref{eq:1-3})--(\ref{eq:2-3}), we should notice that both the coefficients
 $\mathbb{A}^{ij}$ in operators,
 and the time derivative coefficients, can be discontinuous.
 For instance, only if the capillary pressure satisfies $p_c'(S_g=0) =
+\infty$, as in the van Genuchten  model, then all the coefficients
in (\ref{eq:1-3})--(\ref{eq:2-3}) are continuous;  but if this
condition is not satisfied they will be discontinuous.
 \end {remark}

An other variant for the replacement of the saturation by
$\rho_{tot}^h$,  is presented in \cite{TA02-MP01},where relation
 (\ref{def-RHO})is written
\begin{equation}
    \rho_{tot}^{h}= (1-S_g) C_h p_g +S_g C_v p_g,
    \label{panf:1}
\end{equation}
and is then extended to both the two-phase and the one-phase region
by making $p_g=p_l$ in the liquid saturated region (without the gas
phase). This is leading, in the one-phase region where $S_g=0$ and
the Henry law does not apply, to extend   the  gas saturation by
negative values (still defined by equation (\ref{panf:1}) as a
function of the pressure and the total hydrogen concentration).
After a necessary and ad hoc  extension of the  permeability and
capillary pressure curves, out of the usual positive values of
saturation, it is then possible to modeling both the one-phase flow
and the
 two-phase flow with the same system of equations written with this extended saturation as  main
 unknown,
while  using actually the total hydrogen concentration
$\rho_{tot}^{h}$.

\subsection{Modeling based on the hydrogen concentration in the liquid phase, $\rho_l^h$ }

We have seen that the variables $p_l$ and $\rho_{tot}^{h}$,
introduced in the last section, can describe simply the flow system,
both in the one-phase and in the two-phase regions, independently of
the presence of diffusion or capillary forces. But if  we assume
that the capillary forces are present we can choose an other change
of variables in order to have a system of equations with continuous
coefficients. Namely, using the inverse of the capillary pressure
function, we may define the phase saturation as function of the
hydrogen mass concentration in the liquid, $\rho_l^h$, and of the
liquid pressure, $p_l$; and hence use them  as main unknowns. With
these two variables, $\rho_l^h$ and $p_l$  the two systems
(\ref{consSl})--(\ref{eq:5}) and (\ref{un:eq:12})--(\ref{un:eq:34})
are  transformed in a single system of equations able to describe
both liquid saturated and unsaturated flow.
\begin{figure}
\centering
\includegraphics[angle=0, height=1.6in]{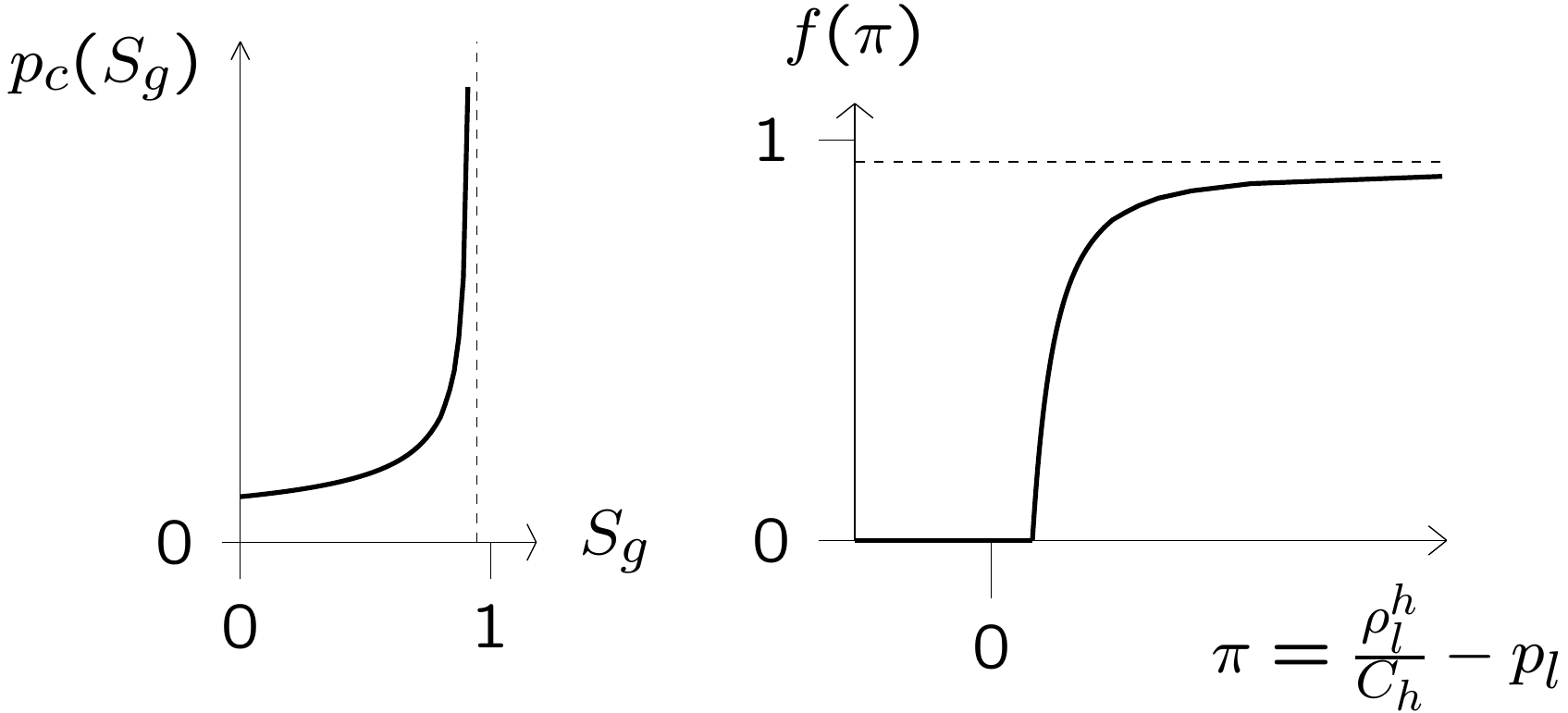}
\caption{Capillary pressure curve, $p_c =p_g - p_l$,  and inverse
function } \label{P-C}
\end{figure}

Since the capillary pressure curve  $S_g\mapsto p_c(S_g)$ is a
strictly increasing function  we can define an inverse function
$f:\mathbb{R}\rightarrow[0,1]$, (see Fig.  \ref{P-C}), by
\begin{equation} \label{eq:II:def_f}
    f(\pi) = \left\{\begin{array}{ll}
        p_c^{-1}(\pi) & \text{if } \pi\geq p_c(0)\\
        0 &\text{otherwise.}
    \end{array}\right.
\end{equation}
By definition of the function  $f $, using (\ref{Henry}) and
(\ref{seuil}),  we have:
\begin{equation} \label{eq:II:Sg=f}
    f\left(\frac{\rho_l^h}{C_h}-p_l\right)= S_g,
\end{equation}
and it is then possible to compute  the gas saturations, $S_g$, from
$p_l$ and $ \rho_l^h$. These two variables being  well defined in
both the one and two-phase regimes, we  will now use them as
principal unknowns.

Equations~\eqref{consSl}-\eqref{eq:5} with  unknowns $p_l$ and $\rho_l^h$ can be written as:
\begin{align}
    -\Phi\rho_w^{std}\der{}{t}\left( f\left(\frac{\rho_l^h}{C_h}-p_l\right)\right)&-
        \dv\big( \Ab^{1,1}\gr p_l + \Ab^{1,2}\gr \rho_l^h + B_1\Kb\gb \big)
        = {\mathcal{F}^w} \label{eq:II:sys_p-r_1}\\
     \Phi\der{}{t}\big( a^{*}\circ f\left(\frac{\rho_l^h}{C_h}-p_l\right)\rho_l^h \big)&-
        \dv\big( \Ab^{2,1}\gr p_l + \Ab^{2,2}\gr \rho_l^h + B_2 \Kb\gb\big)
        = {\mathcal{F}^h} \label{eq:II:sys_p-r_2}
\end{align}
where the coefficients are given by the following formulas:
\begin{align}
\Ab^{1,1} &= \lambda_l(S_g)\rho_w^{std}\Kb,\quad \Ab^{1,2} = -{\Phi (1-S_g)} D \Ib\label{eq:II:A11}\\
\Ab^{2,1} &= \lambda_l(S_g) \rho_l^h \Kb ,\quad
\Ab^{2,2} = \lambda_g(S_g) \frac{C_v}{C_h^2}\rho_l^h \Kb
         + {\Phi (1-S_g)} D \Ib \label{eq:II:A22}\\ 
B_1     &= - \lambda_l(S_g)\rho_w^{std}(\rho_w^{std}+\rho_l^h) \label{eq:II:B1}\\
B_2     &= -\lambda_l(S_g)\rho_l^h(\rho_w^{std}+\rho_l^h)
         - \lambda_g(S_g) \frac{C_v^2}{C_h^2}(\rho_l^h)^2\ \label{eq:II:B2}
\end{align}
with
\begin{equation} \label{eq:II:compl_sys}
    a^*(S_g)= \frac{a(S_g)}{C_h} = 1+(\frac{C_v}{C_h}-1)S_g.
\end{equation}
If we consider first, equation~\eqref{eq:II:sys_p-r_2}, we may write
it  as
\begin{align*}
    \Phi \left(a^{*}(S_g)+\rho_l^h\der{a^{*}(S_g)}{\rho_l^h}\right)\der{\rho_l^h}{t}
    -\dv\big( \Ab^{2,1}\gr p_l &+  \Ab^{2,2}\gr \rho_l^h + B_2\Kb\gb \big)\\
    &+\Phi \rho_l^h\der{a^{*}(S_g)}{p_l}\der{p_l}{t}
    = {\mathcal{F}^h}~.
    \end{align*}
Moreover, from \eqref{eq:II:compl_sys} and because $f$ and $f'$ are
positive, we have
\begin{align*}
    a^{*}(S_g)+\rho_l^h\der{a^{*}(S_g)}{\rho_l^h}
    &=1+(\frac{C_v}{C_h}-1) \left( f\left(\frac{\rho_l^h}{C_h}-p_l\right)
    + \frac{\rho_l^h}{C_h}f'\left(\frac{\rho_l^h}{C_h}-p_l\right)
    \right)
    \geq1 ;
\end{align*}

and if the diffusion is not neglected, we have definite positiveness
of the quadratic form $\Ab^{2,2}~$, in
equation~\eqref{eq:II:sys_p-r_2}; i.e. for any $\boldsymbol{\xi}\neq
0$,
$$
( \Ab^{2,2}\boldsymbol{\xi}\cdot\boldsymbol{\xi}) =
\lambda_g(S_g)\frac{C_v}{C_h^{2}}\rho_l^h
\Kb\boldsymbol{\xi}\cdot\boldsymbol{\xi} + {\Phi (1-S_g)} D
|\boldsymbol{\xi}|^2  >  0,
$$
and therefore equation~\eqref{eq:II:sys_p-r_2} is strictly parabolic
in $\rho_l^h$.\\
If we develop, equation~\eqref{eq:II:sys_p-r_1}  as follows:
\begin{align*}
\Phi\rho_w^{std} f'\left(\frac{\rho_l^h}{C_h}-p_l\right) \der{p_l}{t}
-\dv\big( \Ab^{1,1}\gr p_l &+ \Ab^{1,2}\gr \rho_l^h + B_1\Kb\gb \big)\\
&-\frac{\rho_w^{std}}{C_h}\Phi
f'\left(\frac{\rho_l^h}{C_h}-p_l\right) \der{\rho_l^h}{t} =
{\mathcal{F}^w};\ \label{eq:II:sys_p-r_10}
\end{align*}
we have, for any $\boldsymbol{\xi}$,
$$
\lambda_l(S_g)\rho_w^{std}\Kb\boldsymbol{\xi}\cdot\boldsymbol{\xi}
\geq 0,
$$
and then positiveness of
$(\Ab^{1,1}\boldsymbol{\xi}\cdot\boldsymbol{\xi})$ and of
$(\Ab^{2,1}\boldsymbol{\xi}\cdot\boldsymbol{\xi})$.

Moreover,
$$\Phi\rho_w^{std}
f'\left(\frac{\rho_l^h}{C_h}-p_l\right)\geq0 .$$

However, equations in system
\eqref{eq:II:sys_p-r_1}-\eqref{eq:II:sys_p-r_2} are not uniformly
parabolic/elliptic for the pressure $p_l$, because the coefficients,
$\Ab^{1,1},~\Ab^{2,1} $, in front of $\nabla p_l$ in
\eqref{eq:II:sys_p-r_1}-- \eqref{eq:II:sys_p-r_2} tend to zero as
$S_g\to 1$.

\begin {remark}
It is worth noticing that this system
\eqref{eq:II:sys_p-r_1}-\eqref{eq:II:sys_p-r_2}, with variables
$p_l$ and $ \rho_l^h$, has interesting properties for numerical
simulations in strongly heterogeneous porous media. These two
variables are continuous through interfaces separating different
porous media with different rock types (different absolute
permeability, different capillary  and  permeability curves), as we
will see  in \ref{Sec: 4.3}; which is absolutely not the case for
the variables $p_l$ and $ \rho_{tot}^h$. An other advantage is the
continuity of all the coefficients $\mathbb{A}^{i,j}$, in
\eqref{eq:II:sys_p-r_1}--\eqref{eq:II:sys_p-r_2} and the continuity
of $f$ in \eqref{eq:II:sys_p-r_2}
, even if $p_c'(S_g=0) = +\infty $.

\end {remark}


\section{Numerical experiments}
\label{Sec: 4}
 In this last section, we present four numerical
tests specially designed for illustrating the ability  of the model
described by equations
(\ref{eq:II:sys_p-r_1})-(\ref{eq:II:sys_p-r_2}) to deal with gas
phase appearance and disappearance. Although all the computations
were done using the variables, $p_l$ and $\rho_l^h$, we are also
displaying, for each test,  the Saturation and Pressure level
curves. These two last quantities are obtained after a post
processing step using the Capillary Pressure law (\ref{capillary}),
equations (\ref{eq:II:Sg=f}), Henry's law
(\ref{Henry}), and following the constraints
(\ref{capillaryconstraint}) and (\ref{seuil}) (see Figure \ref{P-D}).\\

\begin{figure}%
\includegraphics[width=.7\columnwidth]{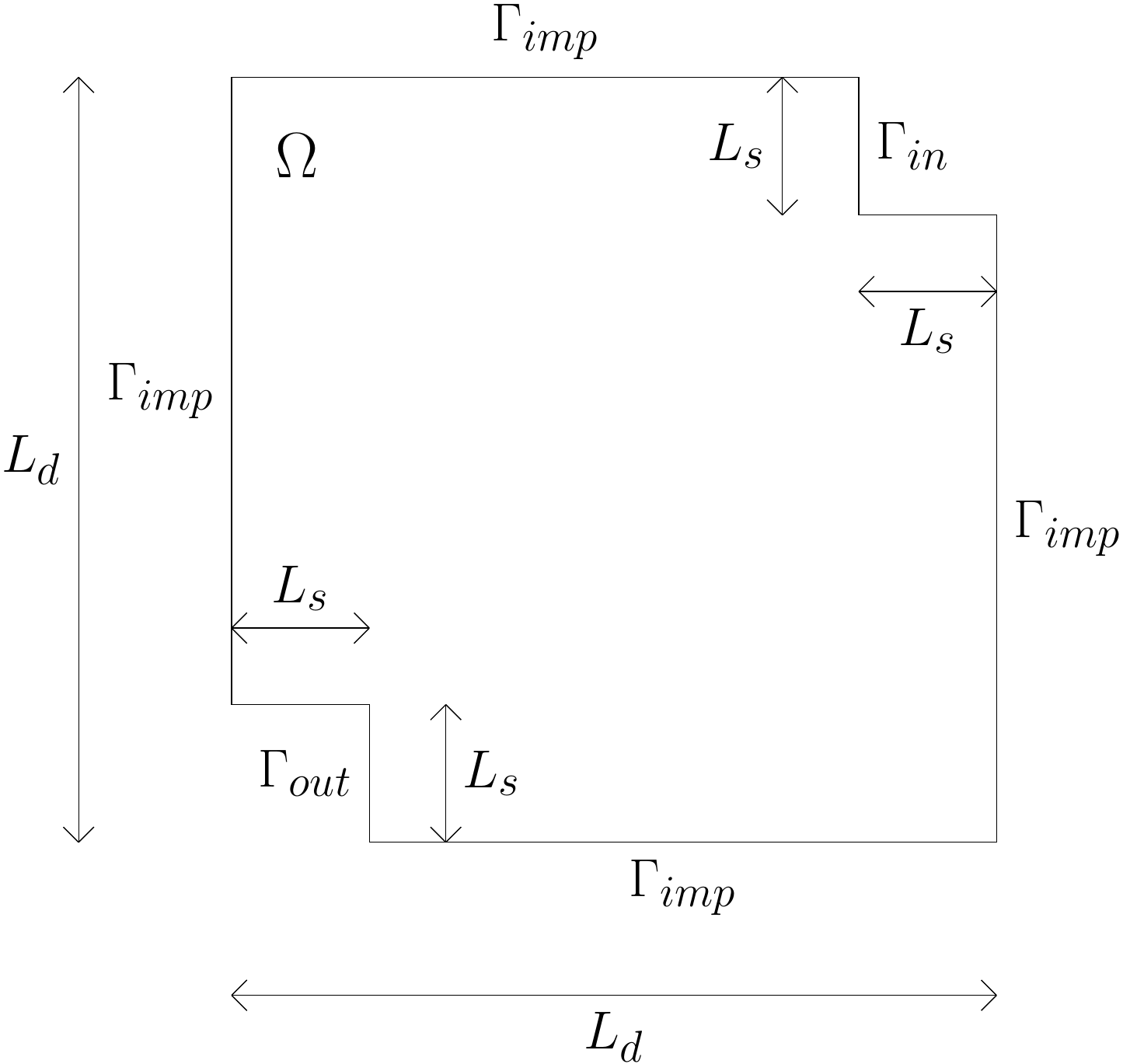}%
\caption{Test case number 1: Geometry a the 2-D porous domain, $\Omega$.}%
\label{fig:geom2D}%
\end{figure}

 The
first test focuses on the gas phase appearance produced by injecting
pure hydrogen in a 2-D homogeneous porous domain $\Omega$ (see
Figure~\ref{fig:geom2D}), which is initially liquid saturated by
pure water(water saturated).

\begin{figure}%
\includegraphics[width=.7\columnwidth]{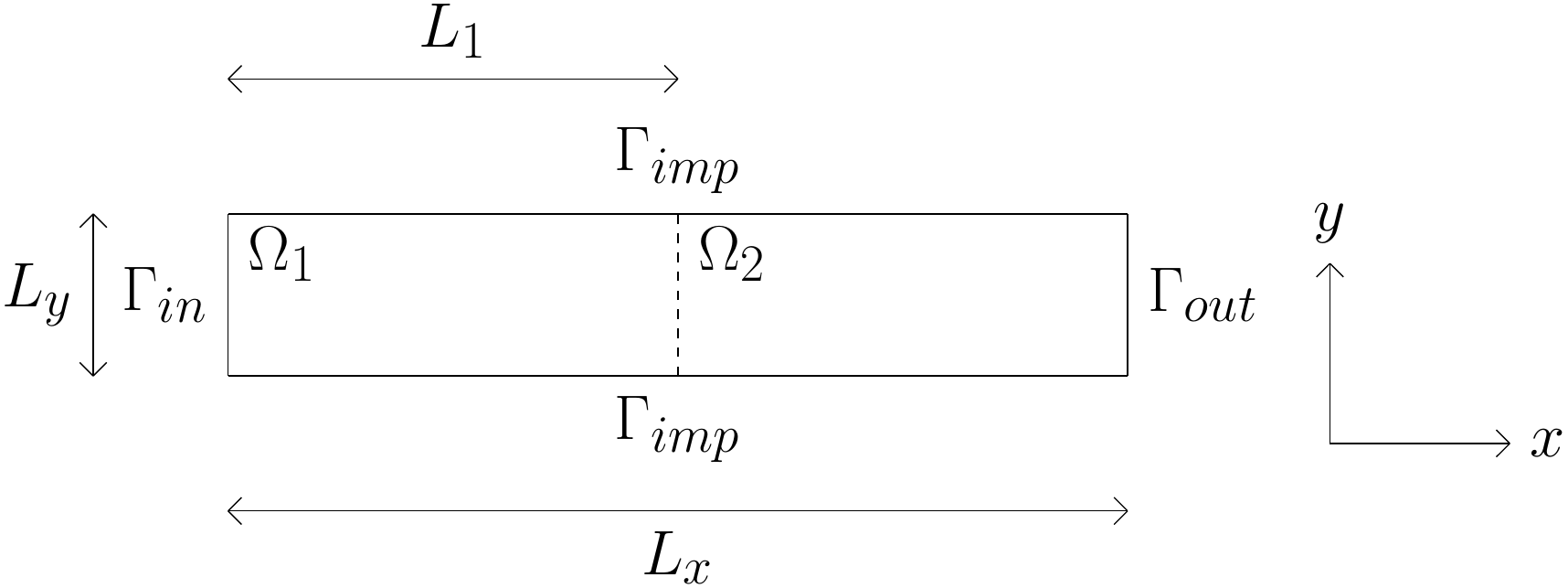}%
\caption{Test cases number 2, 3 and 4: Geometry of the quasi-1D porous domain, $\overline{\Omega} = \overline{\Omega_1} \cup \overline{\Omega_2}$.}%
\label{fig:geom1D}%
\end{figure}

 Because the main goal of all these numerical experiments is to test  the model  efficiency, for
 describing the phase appearance or disappearance, the porous domain geometry does not really matter and we will use  a porous domain with a simple
 geometry. Consequently, we choose a simple, quasi-1D, porous domain (see Figure~\ref{fig:geom1D}) for the all  three next tests .

 The test case number 2 is more complex, it shows local
disappearance of the gas phase created by injecting  pure hydrogen
in a homogeneous unsaturated porous medium (initially both phases,
liquid and gas, are present everywhere).

 The two last tests  aim is to focus on the
main challenges in simulating the flow crossing the engineered
barriers, located around the waste packages.  In the test case
number 3, the porous medium domain is split in two parts with
different and highly contrasted rock types, and like in the first
one, the gas phase appearance is produced by injecting pure hydrogen
in an initially water saturated porous domain. The test case number
4 addresses the evolution of the phases, from an initial phase
disequilibrium to a stabilized stationary state, in a closed porous
domain (no flux boundary conditions).
\begin{table}[htb] \centering
    \begin{tabular}{|c|rc|}
        \hline
        Parameter & \multicolumn{2}{c|}{Value} \\
        \hline
        $\theta$ & $303$ & $K$ \\
        $D_l^h$ & $3\;10^{-9}$ & $m^2/s$ \\
        $\mu_l$ & $1\;10^{-3}$ & $Pa.s$ \\
        $\mu_g$ & $9\;10^{-6}$ & $Pa.s$ \\
        $H(\theta=303K)$ & $7.65\;10^{-6}$ & $mol/Pa/m^3$ \\
        $M^w$ & $10^{-2}$ & $kg/mol$ \\
        $M^h$ & $2\;10^{-3}$ & $kg/mol$ \\
        $\rho_w^{std}$ & $10^3$ & $kg/m^3$ \\
        \hline
    \end{tabular}
    \caption{Fluid parameters: phases and components characteristics.}
    \label{tab:param_wh}
\end{table}
\begin{table}
\center
\begin{tabular}{l|cl|c|}
    & \multicolumn{2}{c|}{ Mesh size range } & Time step range \\
    \hline Test number 1
    & 2 $m$ -- \ 6 $m$ & $^{(*)}$ & $10^2$ years -- \ 5 \ $10^4$ years \\
    \hline Test number 2
    & 1 $m$ & $^{(**)}$ & $10^2$ years -- \ 5 \ $10^3$ years \\
    \hline Test number 3
    & 1 $m$ & $^{(**)}$ & $10^2$ years -- \ 2 \ $10^4$ years \\
    \hline Test number 4
    & 2 \ $10^{-3}$ $m$ & $^{(**)}$ & $0.33$ $s$ -- \   \ $16.7$ $10^3$ $s$ \\
    \hline
    \multicolumn{4}{l}{\small(*) \  Unstructured triangular mesh}\\
    \multicolumn{4}{l}{\small(**)  Regular quadrangular mesh}
\end{tabular}
\caption{Mesh sizes and time steps used  in the different Numerical
Test } \label{tab:NUM}
\end{table}
In all these four test cases, for simplicity,  the porous medium is
assumed to be isotropic, such that $\mathbb{K}=k\mathbb{I}$ with $k$
a positive scalar; and the source terms are assumed to be null:
$\mathcal{F}_w=0$ and $\mathcal{F}_h=0$. As usual in geohydrology,
 the van Genuchten-Mualem model for the capillary pressure law
and the relative permeability functions are used in underground
nuclear waste modeling, i.e. :
\begin{equation} \label{eq:VGM}
    \left\{
    \begin{gathered}
        p_c = P_r\left(S_{le}^{-1/m}-1\right)^{1/n}
        \ ,\ 
        \lambda_{l} = \frac{1}{\mu_l}\sqrt{S_{le}}\left( 1-(1-S_{le}^{1/m})^m \right)^2
        \\\text{and}\quad
        \lambda_{g}=\frac{1}{\mu_g}\sqrt{1-S_{le}}\left( 1-S_{le}^{1/m} \right)^{2m}
        \\\text{with}\quad
        S_{le}=\frac{S_l-S_{l,res}}{1-S_{l,res}-S_{g,res}}
        \quad\text{and}\quad
        m=1-\frac{1}{n}.
    \end{gathered}
    \right.
\end{equation}
Note that in the van Genuchten-Mualem model, we have no entry
pressure, $p_c(0)=0$,
but the presence of an entry pressure will not lead to any difficulty, neither from the mathematical point of view,
 nor for the numerical simulations.
 Concerning the other fluid characteristics, the values of the physical
parameters specific to the phases (liquid and gas) and to the
components (water and hydrogen) are given in
Table~\ref{tab:param_wh}. All the simulations, presented herein,
were performed using the modular code {\it Cast3m}, \cite{castem}.
The differential equations system was first linearized by a
quasi-Newton method and then discretized by a finite volume,
implicit in time, scheme; with the discretization parameters (mesh
size and time step) given in Table~ \ref{tab:NUM}.

\subsection{Numerical Test number 1}
\label{Sec: 4.1}
 The geometry of this test case is given in
Figure~\ref{fig:geom2D}; and the related data  are given in
Table~\ref{tab:ct1}. A constant flux of hydrogen is imposed on the
input boundary,
 $\Gamma_{in}$, while Dirichlet conditions $p_l=p_{l,out}$, $\rho_l^h=0$ are given
 on $\Gamma_{out}$in order to have only the water component on this part of the boundary. The initial conditions, $p_l=p_{l,out}$ and
 $\rho_l^h=0$, are uniform on all the domain, and correspond to a porous domain initially saturated with pure
 water.\\

\begin{figure}[hbtp]
\begin{tabular}{ccc}
    \hspace{4em}$\rho_l^h$\hspace{1ex}
        \includegraphics[height=.1\textheight]{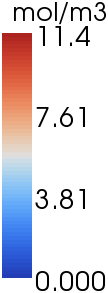}&
    \hspace{4em}$p_l$\hspace{1ex}
        \includegraphics[height=.1\textheight]{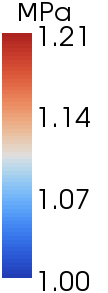}&
    \hspace{4em}$S_g$\hspace{1ex}
        \includegraphics[height=.1\textheight]{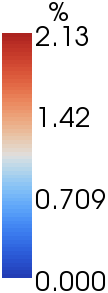}\\
    \includegraphics[width=.3\textwidth]{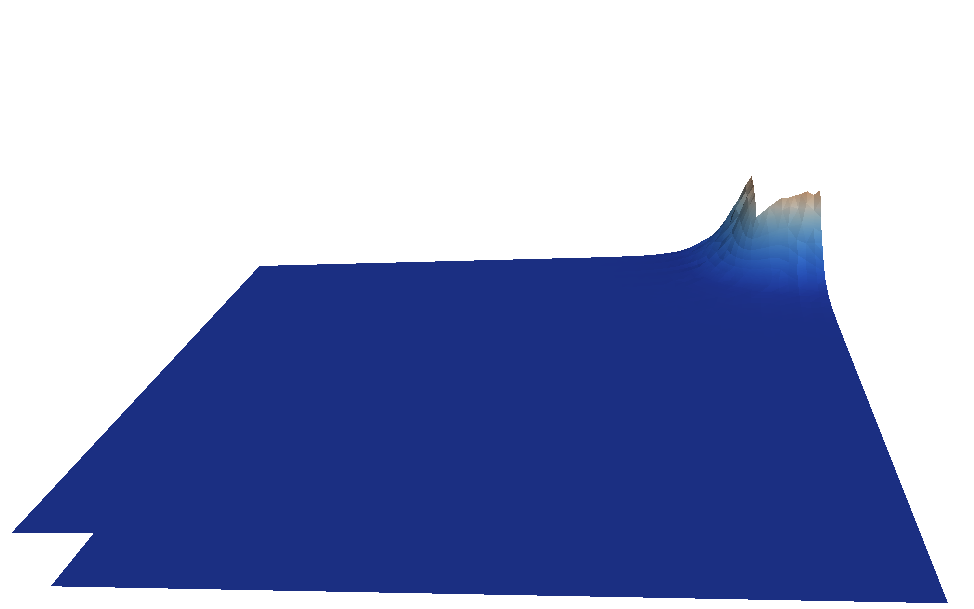}&
    \includegraphics[width=.3\textwidth]{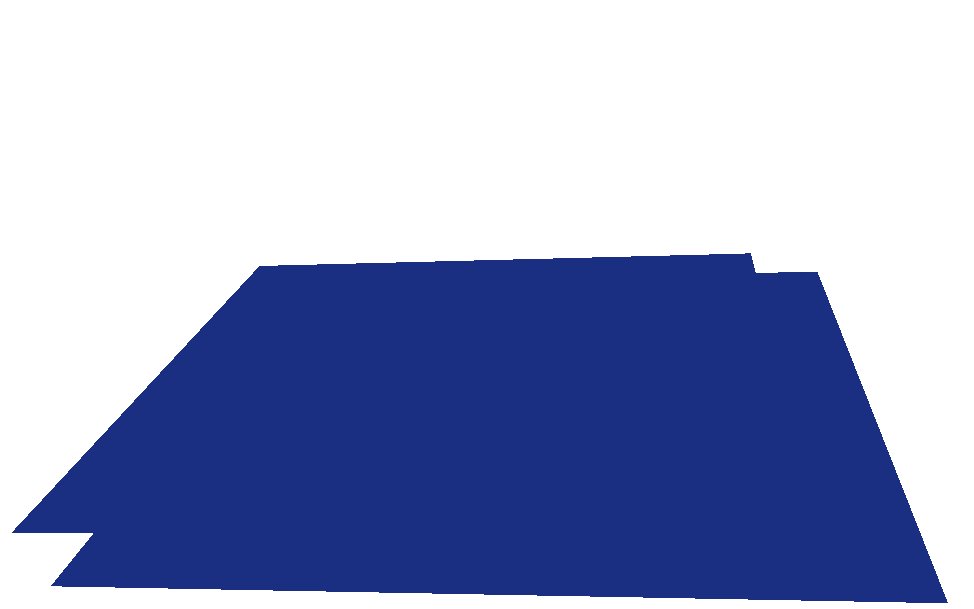}&
    \includegraphics[width=.3\textwidth]{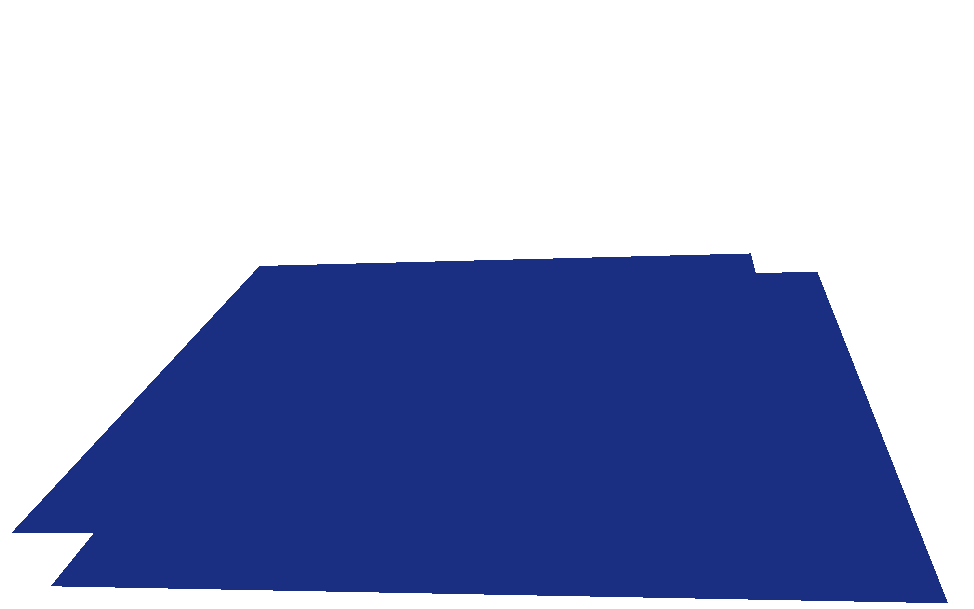}\\
    $\rho_l^h$ at $t=1200$ years&$p_l$ at $t=1200$ years&$S_g$ at $t=1200$ years\\
    \includegraphics[width=.3\textwidth]{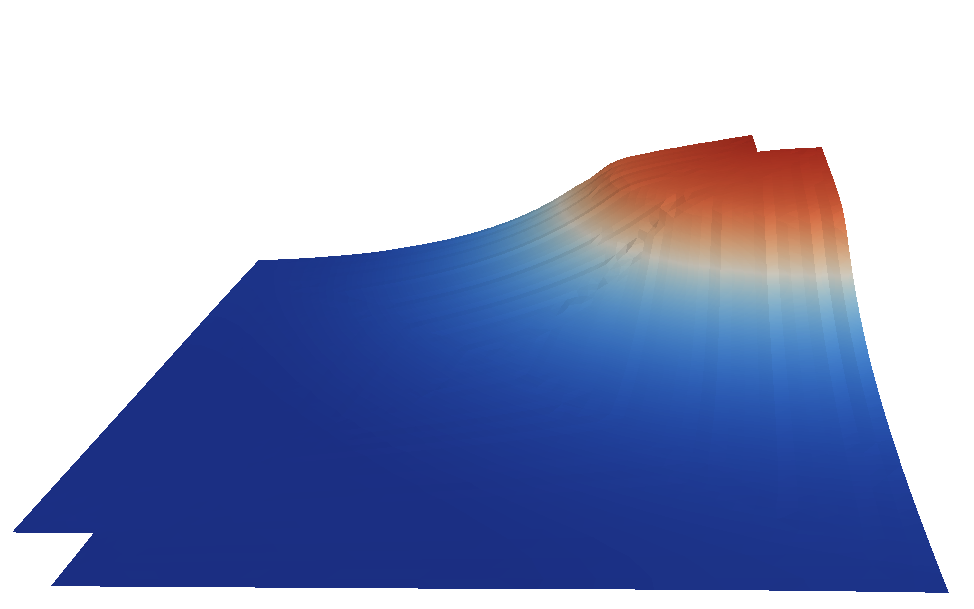}&
    \includegraphics[width=.3\textwidth]{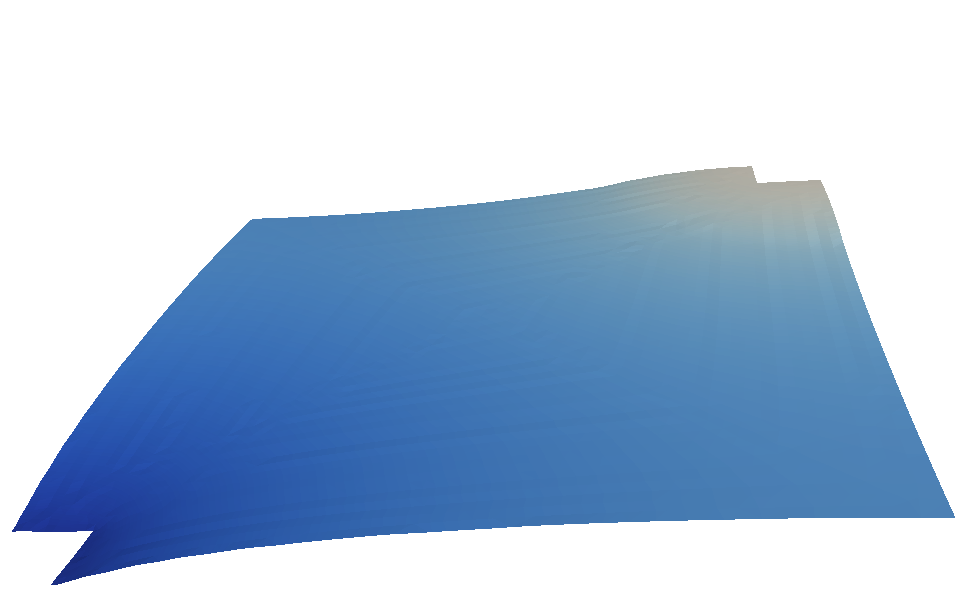}&
    \includegraphics[width=.3\textwidth]{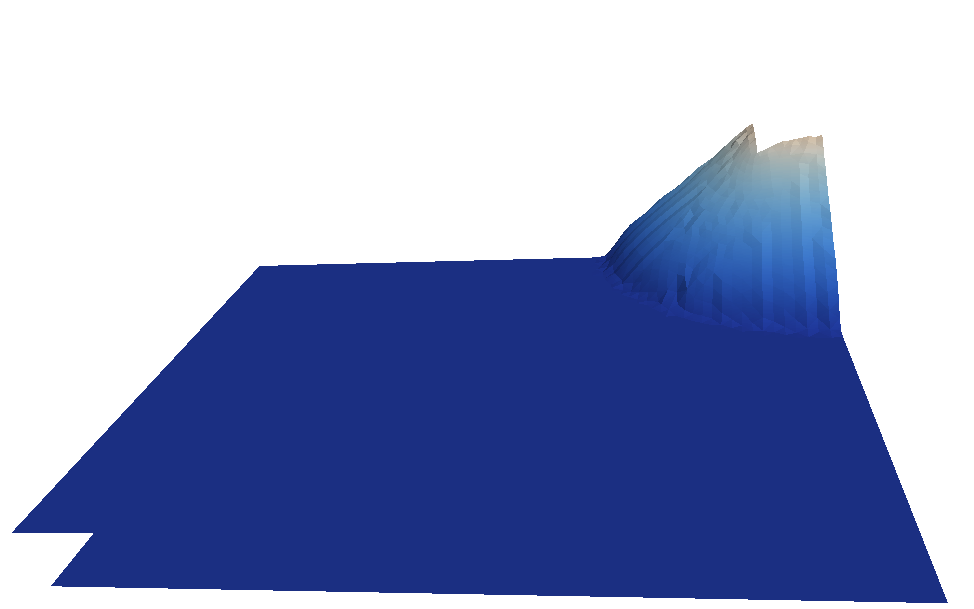}\\
    $\rho_l^h$ at $t=4\:10^4$ years&$p_l$ at $t=4\:10^4$ years&$S_g$ at $t=4\:10^4$ years\\
    \includegraphics[width=.3\textwidth]{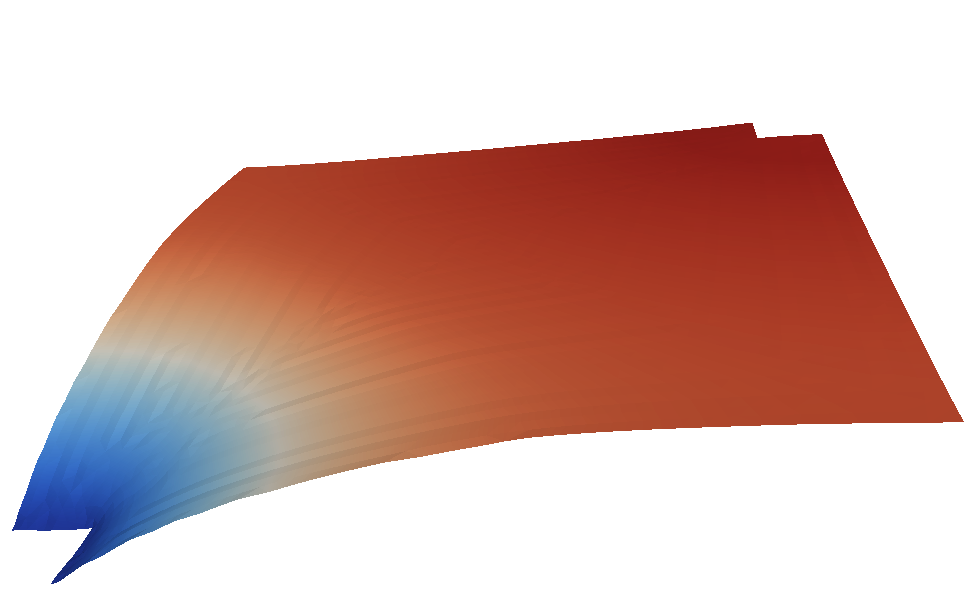}&
    \includegraphics[width=.3\textwidth]{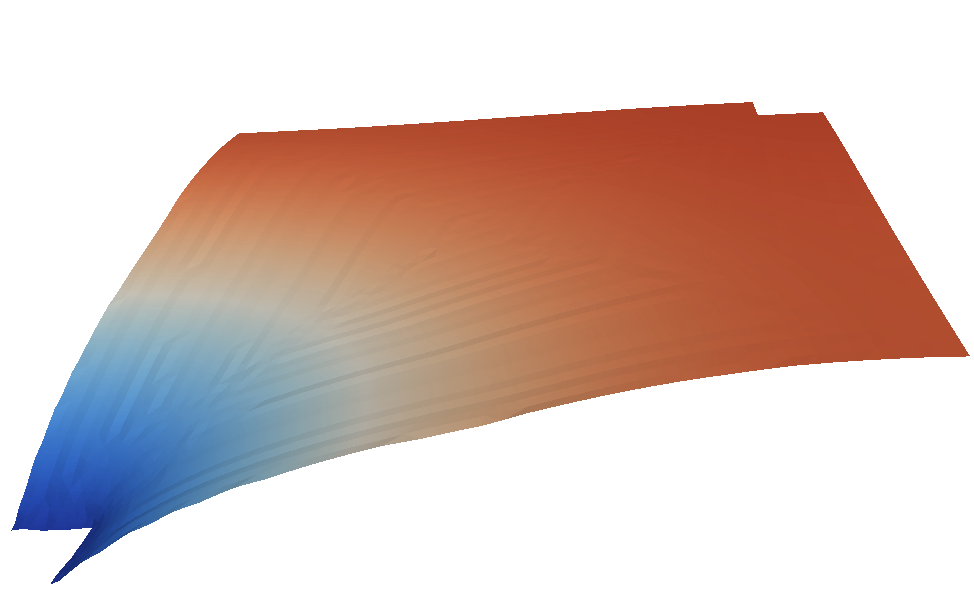}&
    \includegraphics[width=.3\textwidth]{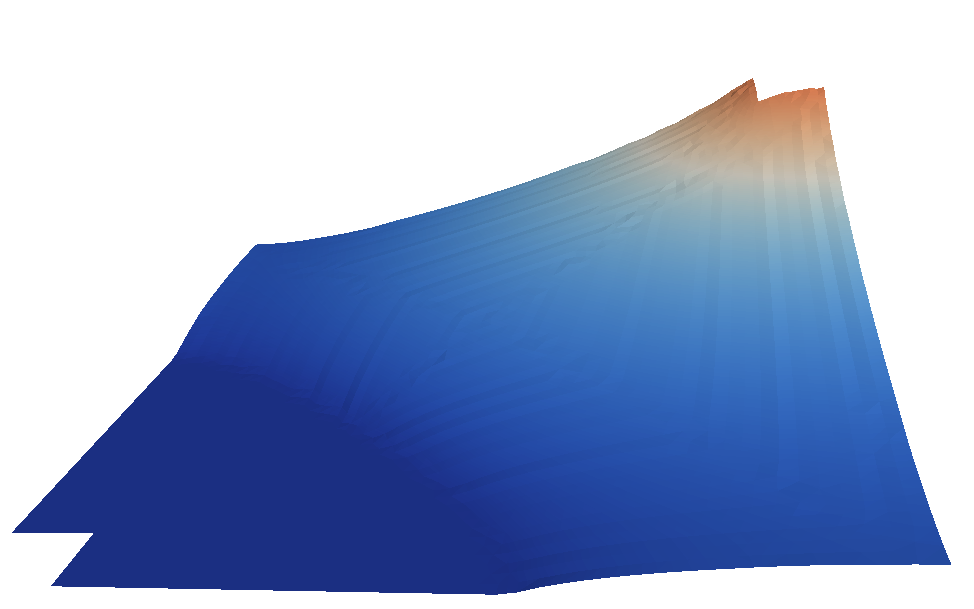}\\
    $\rho_l^h$ at $t=2\:10^5$ years&$p_l$ at $t=2\:10^5$ years&$S_g$ at $t=2\:10^5$ years\\
    \includegraphics[width=.3\textwidth]{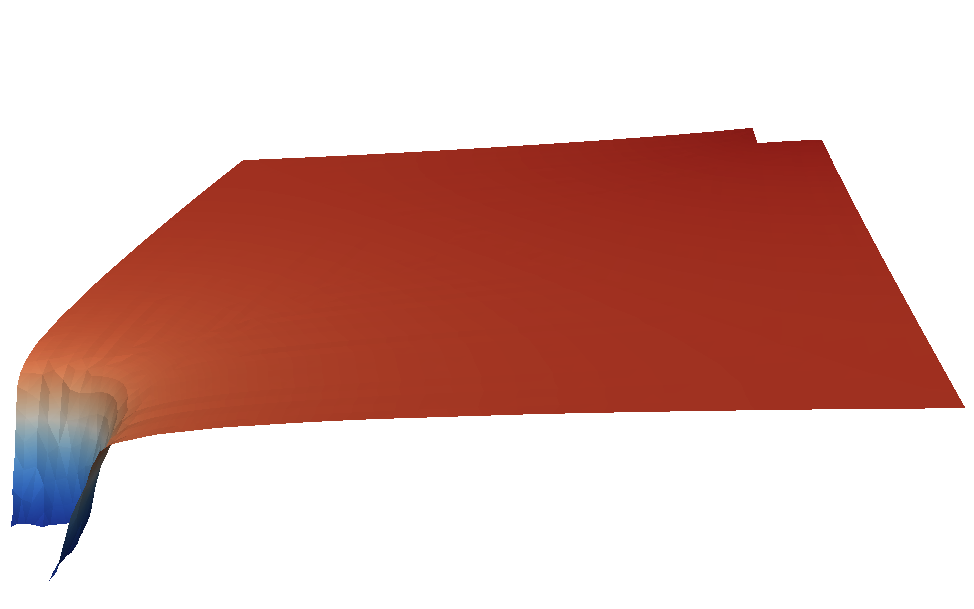}&
    \includegraphics[width=.3\textwidth]{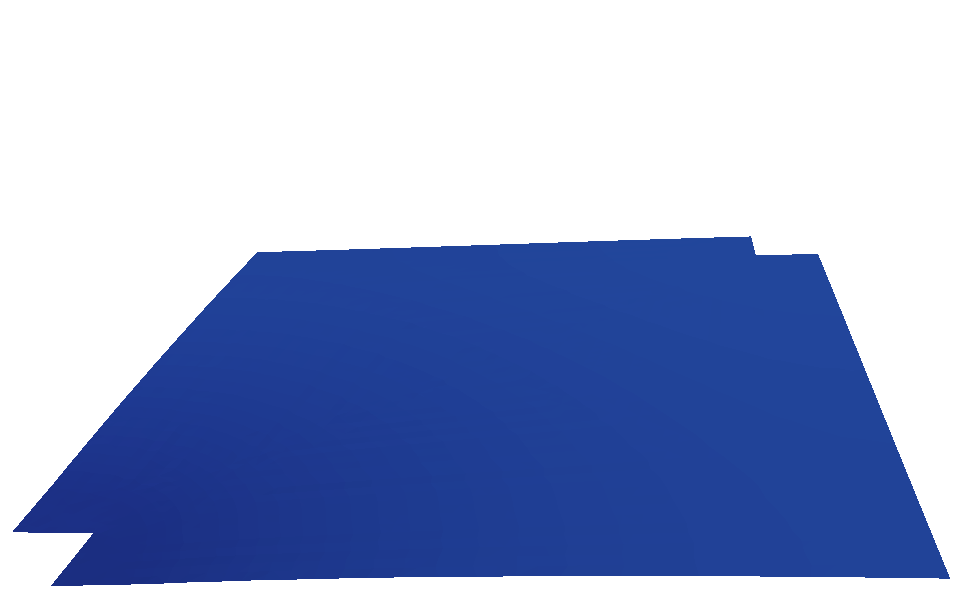}&
    \includegraphics[width=.3\textwidth]{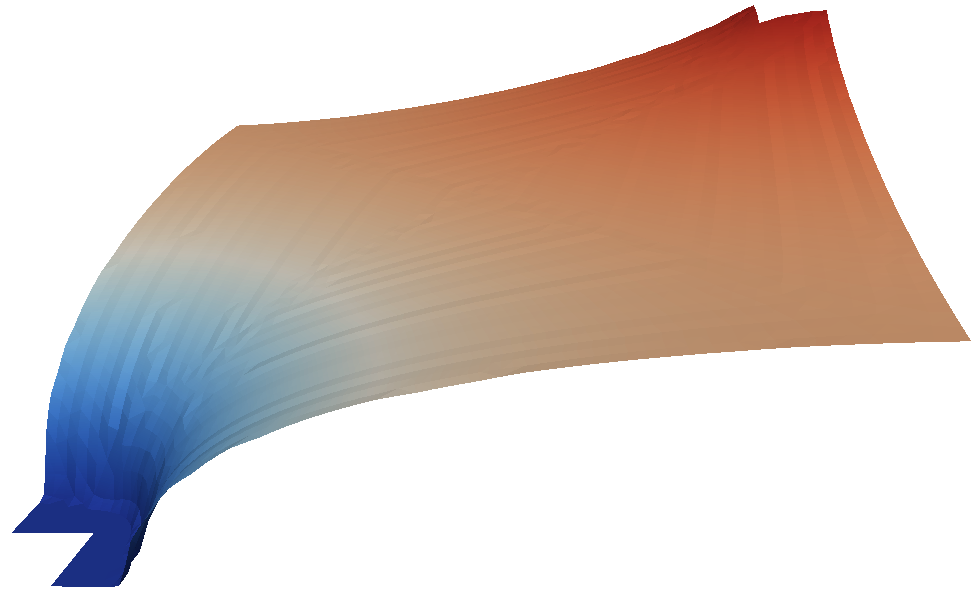}\\
    $\rho_l^h$ at $t=10^6$ years&$p_l$ at $t=10^6$ years&$S_g$ at $t=10^6$ years\\
\end{tabular}
    \caption{Numerical Test case number 1: Evolution of $\rho_l^h$,the hydrogen concentration in the liquid
    phase; $p_l$ the liquid phase pressure; and $S_g$ the gas saturation ; at times
$t=1200, 4\:10^4, 2\:10^5$ and $10^6$ years (from the top to the
bottom). }
    \label{fig:ct1}
\end{figure}

The main steps of the corresponding simulation are presented in
Figure~\ref{fig:ct1}.

\begin{table}[htb] \centering \small
    \begin{tabular}{|c||c|rl||c|rl|}
        \hline
        Boundary conditions & \multicolumn{3}{c||}{Porous medium} & \multicolumn{3}{c|}{Others}\\
        \cline{2-7}
        Initial condition & Param. & \multicolumn{2}{c||}{Value} & Param. & \multicolumn{2}{c|}{Value}\\
        \hline
        $\phi^w\cdot\nu=0$ on $\Gamma_{imp}$ &
        $k$&$5\;10^{-20}$&\hspace{-1.5ex}$m^2$&
        $L_d$&$200$&\hspace{-1.5ex}$m$\\
        $\phi^h\cdot\nu=0$ on $\Gamma_{imp}$ &
        $\Phi$&$0.15$&$\hspace{-1.5ex}(-)$&
        $L_s$&$20$&\hspace{-1.5ex}$m$\\
        $\phi^w\cdot\nu=0$ on $\Gamma_{in}$ &
        $P_r$&$2\;10^6$&\hspace{-1.5ex}$Pa$&
        $p_{l,out}$&$10^6$&\hspace{-1.5ex}$Pa$\\
        $\phi^h\cdot\nu=\mathcal{Q}^h$  on $\Gamma_{in}$ &
        $n$&$1.49$&$\hspace{-1.5ex}(-)$&
        $\mathcal{Q}^h$&$9.28$&\hspace{-1.5ex}$mg/m^2/year$\\
        $p_l=p_{l,out}$ on $\Gamma_{out}$ &
        $S_{l,res}$&$0.4$&$\hspace{-1.5ex}(-)$&
        &&\\
        $\rho_l^h=0$ on $\Gamma_{out}$ &
        $S_{g,res}$&$0$&$\hspace{-1.5ex}(-)$&
        &&\\
        $p_l(t=0)=p_{l,out}$ in $\Omega$ &
        &&&
        &&\\
        $\rho_l^h(t=0)=0$ in $\Omega$ &
        &&&
        &&\\
      \hline
    \end{tabular}
    \caption{Numerical Test case number 1: Boundary and Initial Conditions; porous medium characteristics
     and domain geometry; $\phi^w$ and $\phi^h$ are denoting respectively the water and hydrogen flux.}
      \label{tab:ct1}
\end{table}

We observe in the beginning (see time $t=1200\:$years ~in
Figure~\ref{fig:ct1}) that all the  injected hydrogen through
$\Gamma_{in}$ is totally dissolved in the liquid phase, the gas
saturation stay null on all the domain (there is no gas phase).
During that same period of time: the liquid pressure stay constant,
the liquid phase does not flow, and the hydrogen is transported only
by diffusion of the dissolved hydrogen in the liquid phase.

Later on, the dissolved hydrogen accumulates  around $\Gamma_{in}$
until the dissolved hydrogen concentration $\rho_l^h$ reaches the
threshold $\rho_l^h=C_hp_l$ ( according to Figure~\ref{P-D}and
$p_{c}(0)=0$ in \ref{Rem:1}), at time $t=1600$~ years, when the gas
phase appears in the vicinity of $\Gamma_{in}$. Then this
unsaturated region ( the two-phases, gas and liquid are present
together) progressively expands and the liquid pressure, due to the
compression by the gas phase, increases in the whole porous domain,
causing
 the liquid phase to flow from $\Gamma_{in}$ to $\Gamma_{out}$. Consequently, after
this time,  $t=1600$ years: the  hydrogen is transported by
convection in the gas phase and the dissolved hydrogen is
transported by both convection and diffusion in the liquid phase.
The liquid phase pressure increases globally in the whole domain
until time $t=260\;000$ years (see Figure~\ref{fig:ct1}), and it
starts to decrease in the whole domain until reaching  a uniform and
stationary  state at $t=10^6$ years, in which the water component
flux is null everywhere.

\subsection{Numerical Test number  2}
\label{Sec: 4.2} The geometry and the data of this numerical test
are given in Figure~\ref{fig:geom1D} and Table~\ref{tab:ct2}. The
porous medium is homogeneous and the initial conditions uniform;
there is no need for defining  two parts of the porous domain,
$\Omega_1$ and $\Omega_2$; the parameter $L_1$ will be considered as
null.

In this second test a constant flux of hydrogen is imposed on the
input boundary  $\Gamma_{in}$, while Dirichlet conditions
$p_l=p_{l,out}$, $p_g=p_{g,out}$ are chosen, on $\Gamma_{out}$ ,
such that  $\rho_l^h>C_hp_l$, in order to keep the gas phase (
according to the phase diagram in Figure~\ref{P-D})  present on this
part of the boundary.
 The initial conditions $p_l=p_{l,out}$ and  $\rho_l^h=C_hp_{g,out}$ are uniform and imply
 the presence of the gas phase ( $S_g>0$) in the whole domain. \\

\begin{figure}[hbtp]
    \includegraphics[width=.45\textwidth]{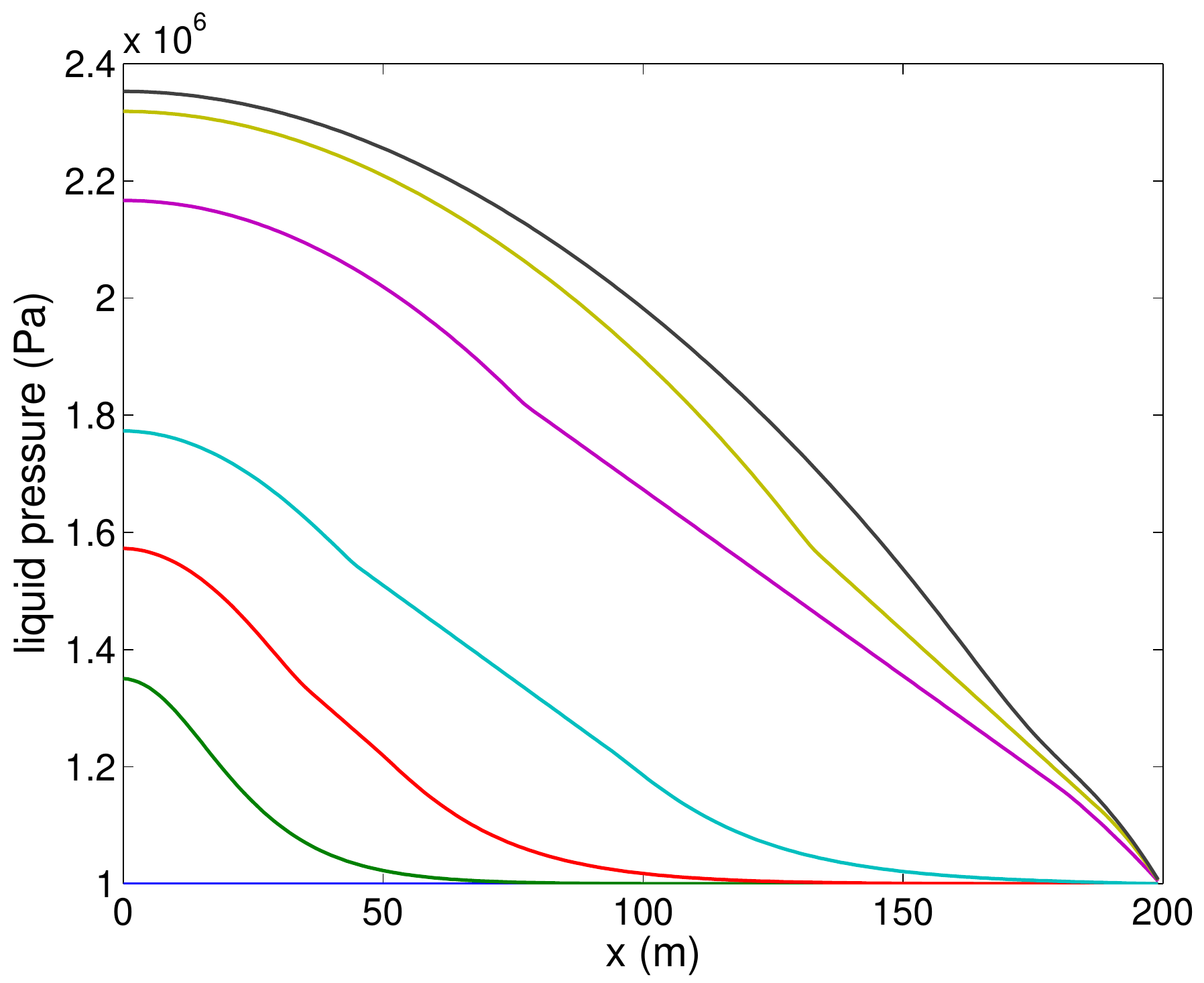}
    \includegraphics[width=.45\textwidth]{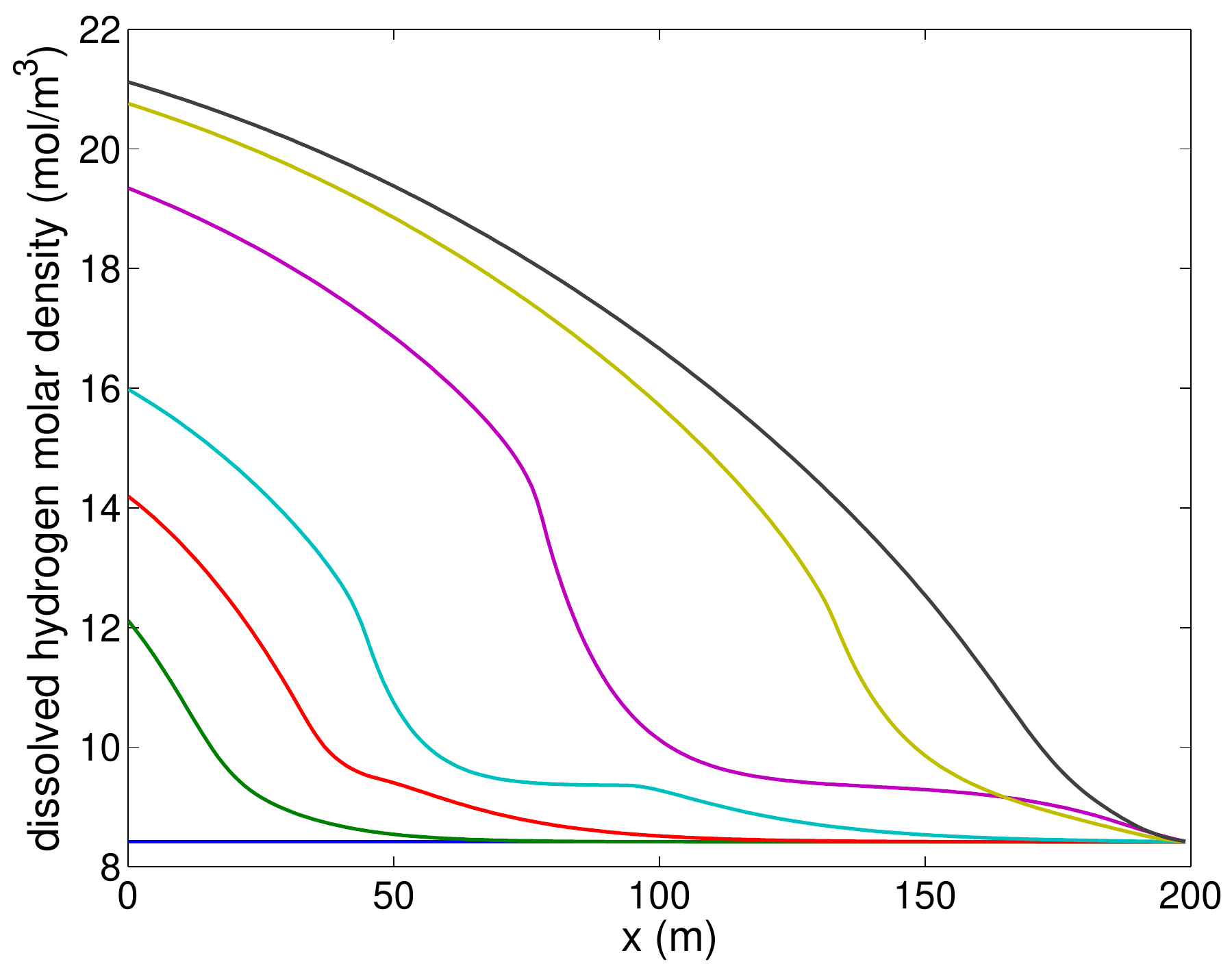}
    \includegraphics[width=.45\textwidth]{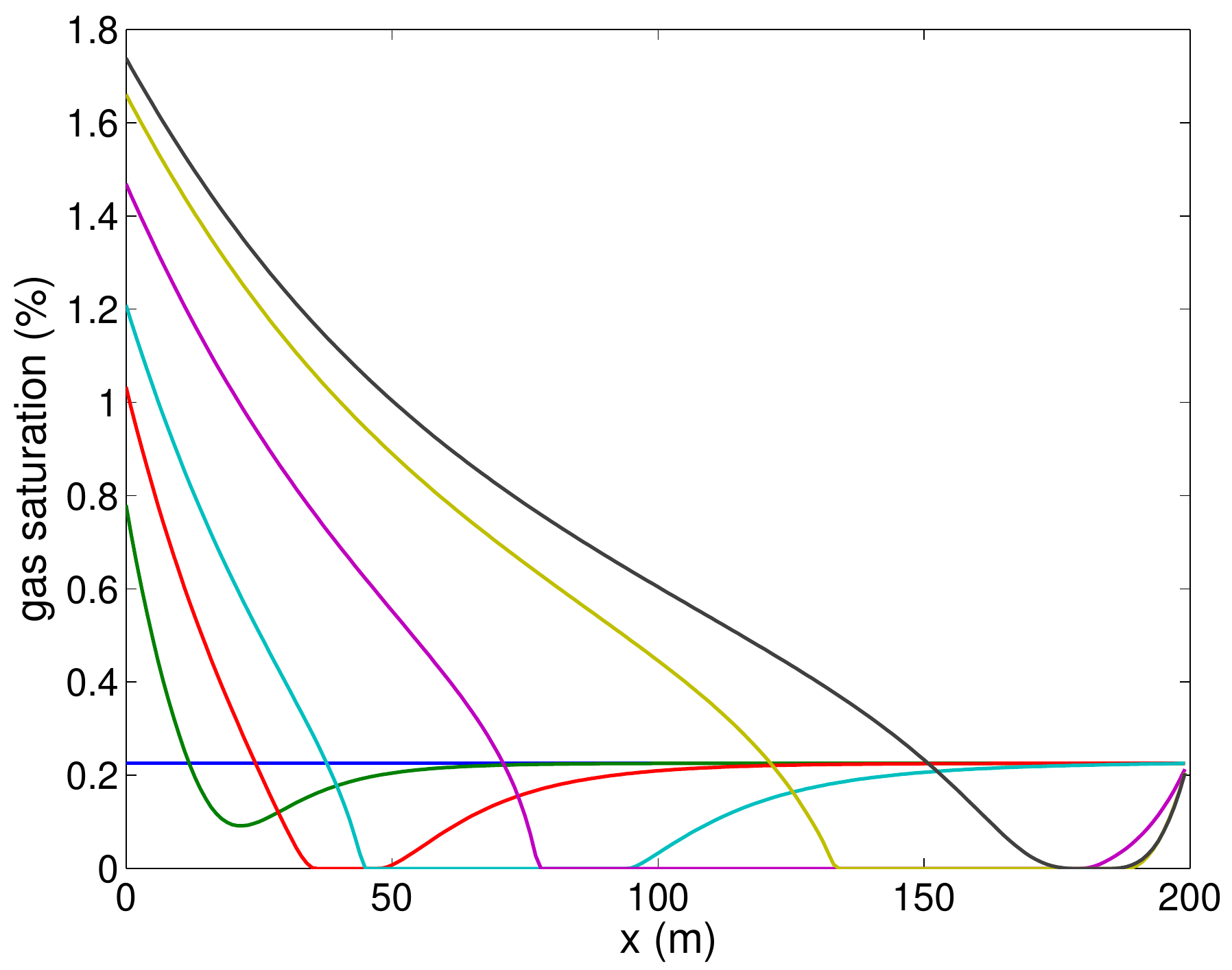}
    \includegraphics[width=.45\textwidth]{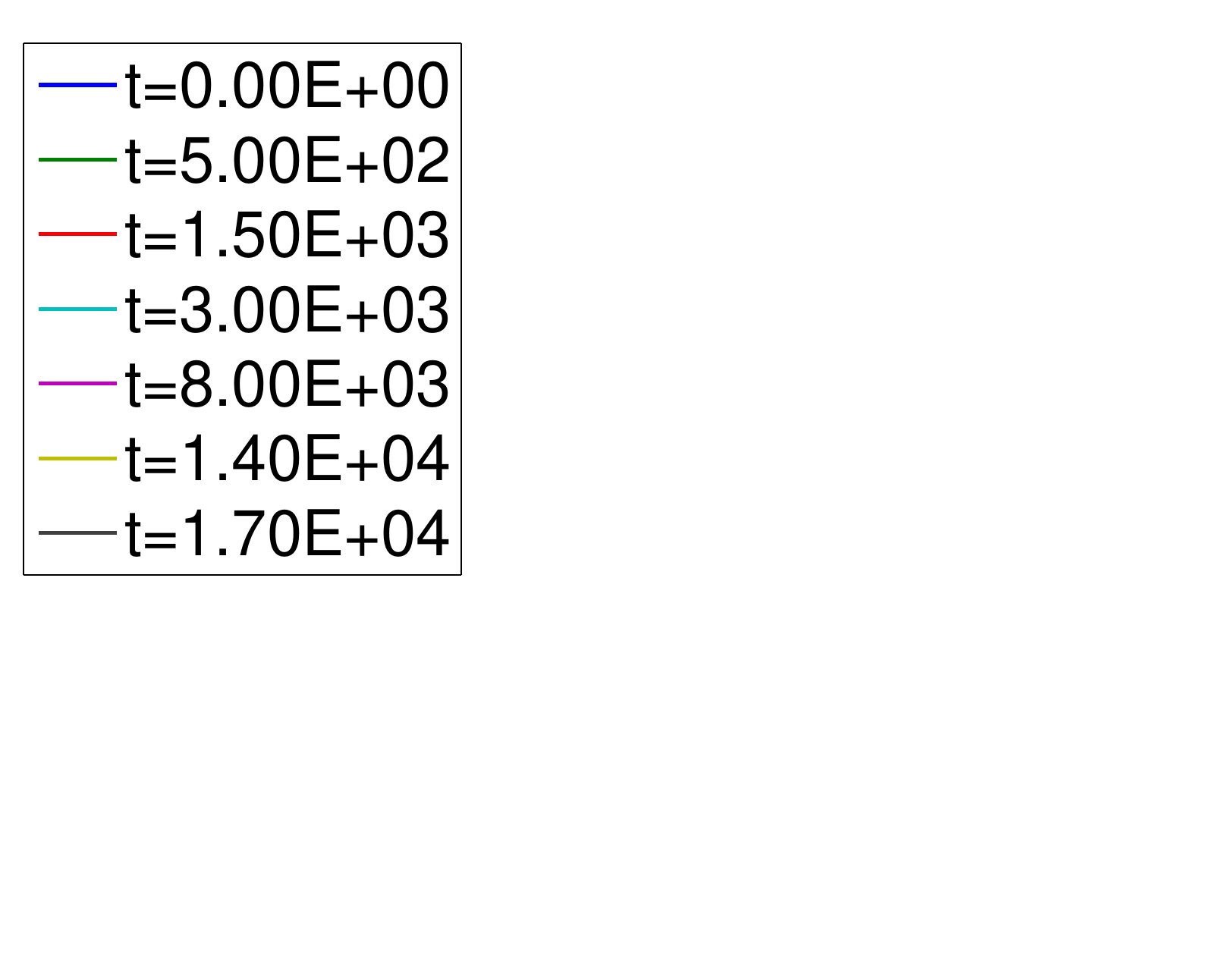}
    \caption{Test case number 2; $L_x= L_2=200$ m:  Time evolution of  $\rho_l^h$ (top right), the dissolved hydrogen
    molar density ($\rho_l^h/M^h$) (top left) and $S_g$ (bottom) profiles ;
    during the first time steps.}
    \label{fig:ct2_1}
\end{figure}

\begin{figure}[hbtp]
    \includegraphics[width=.45\textwidth]{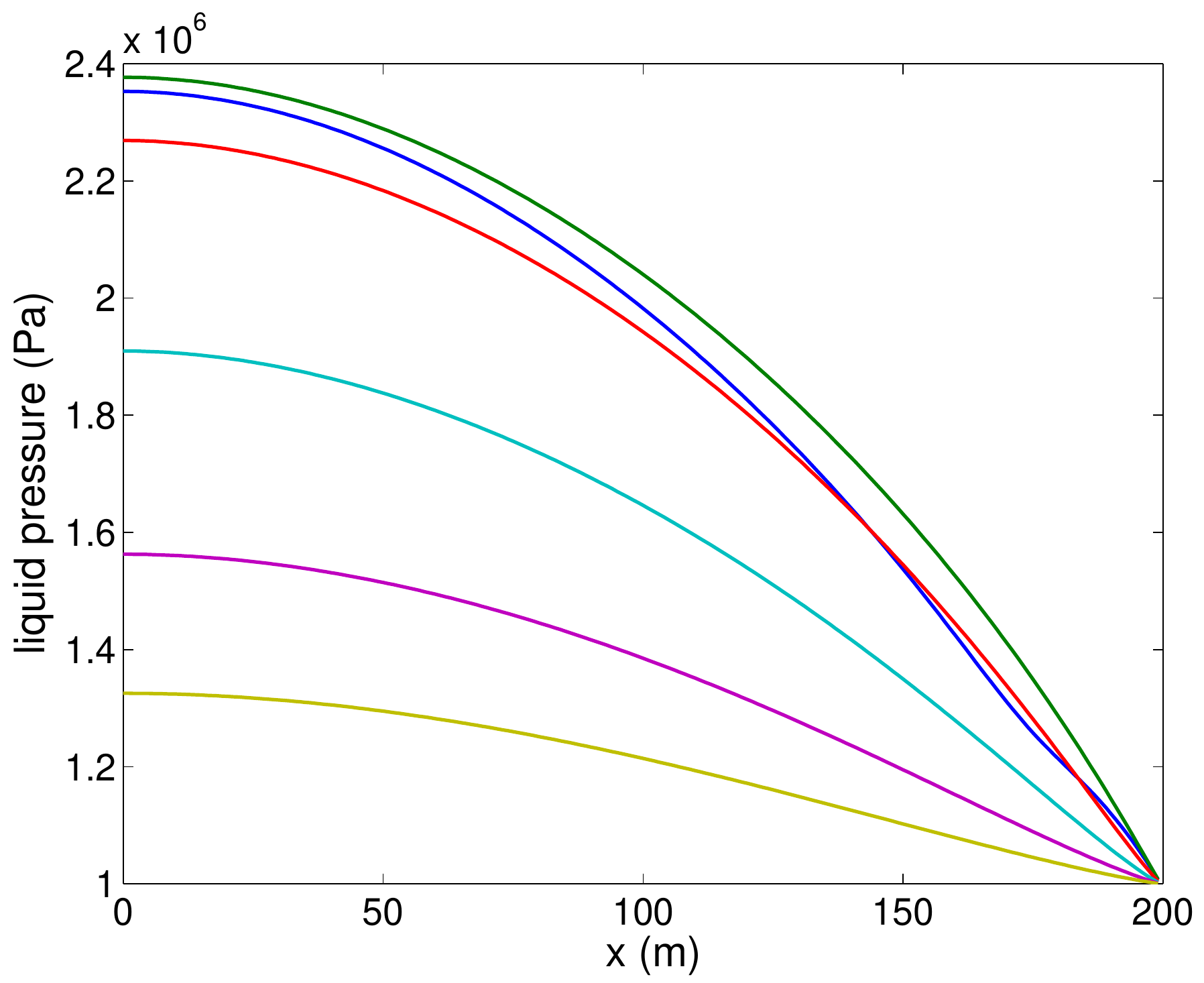}
    \includegraphics[width=.45\textwidth]{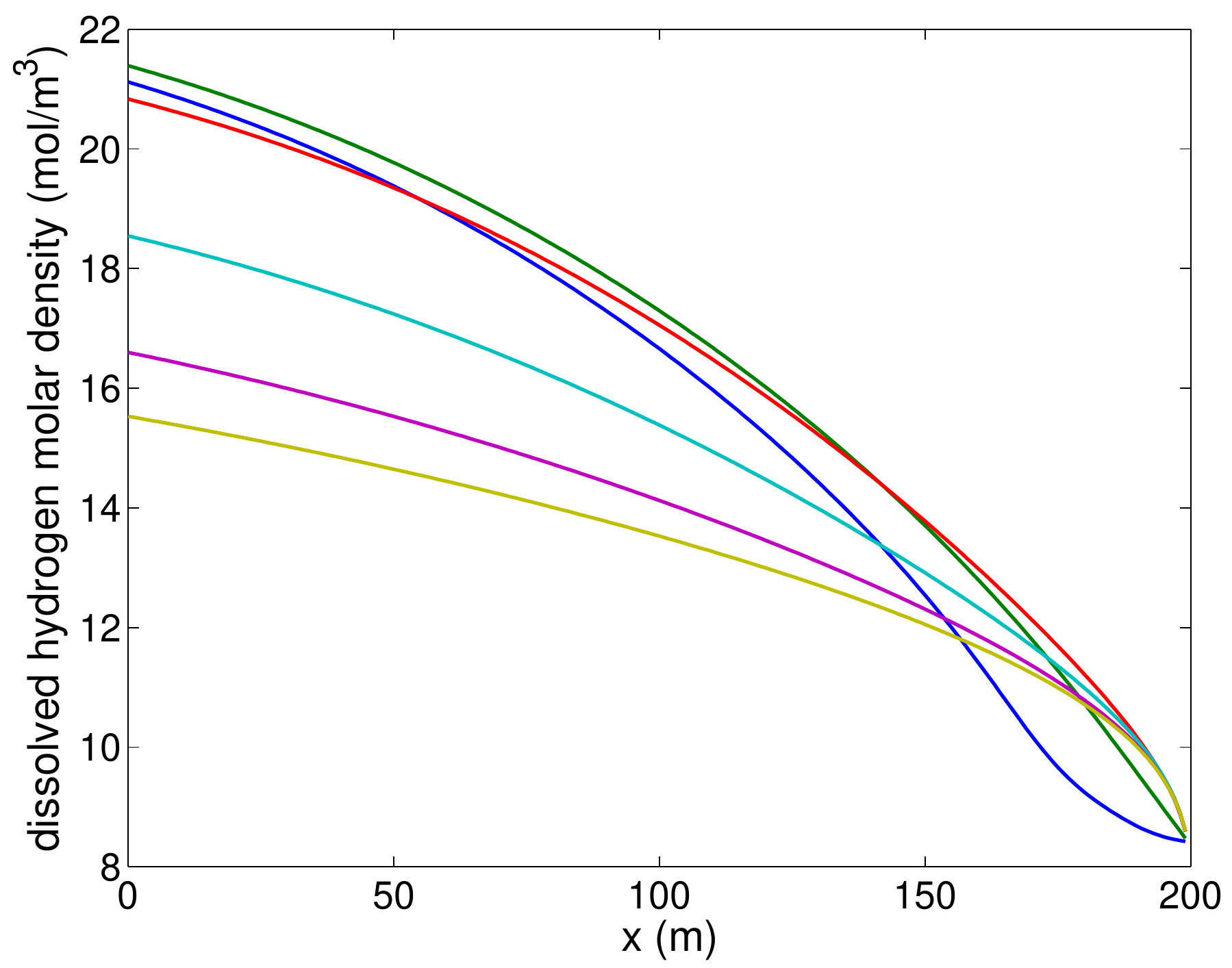}
    \includegraphics[width=.45\textwidth]{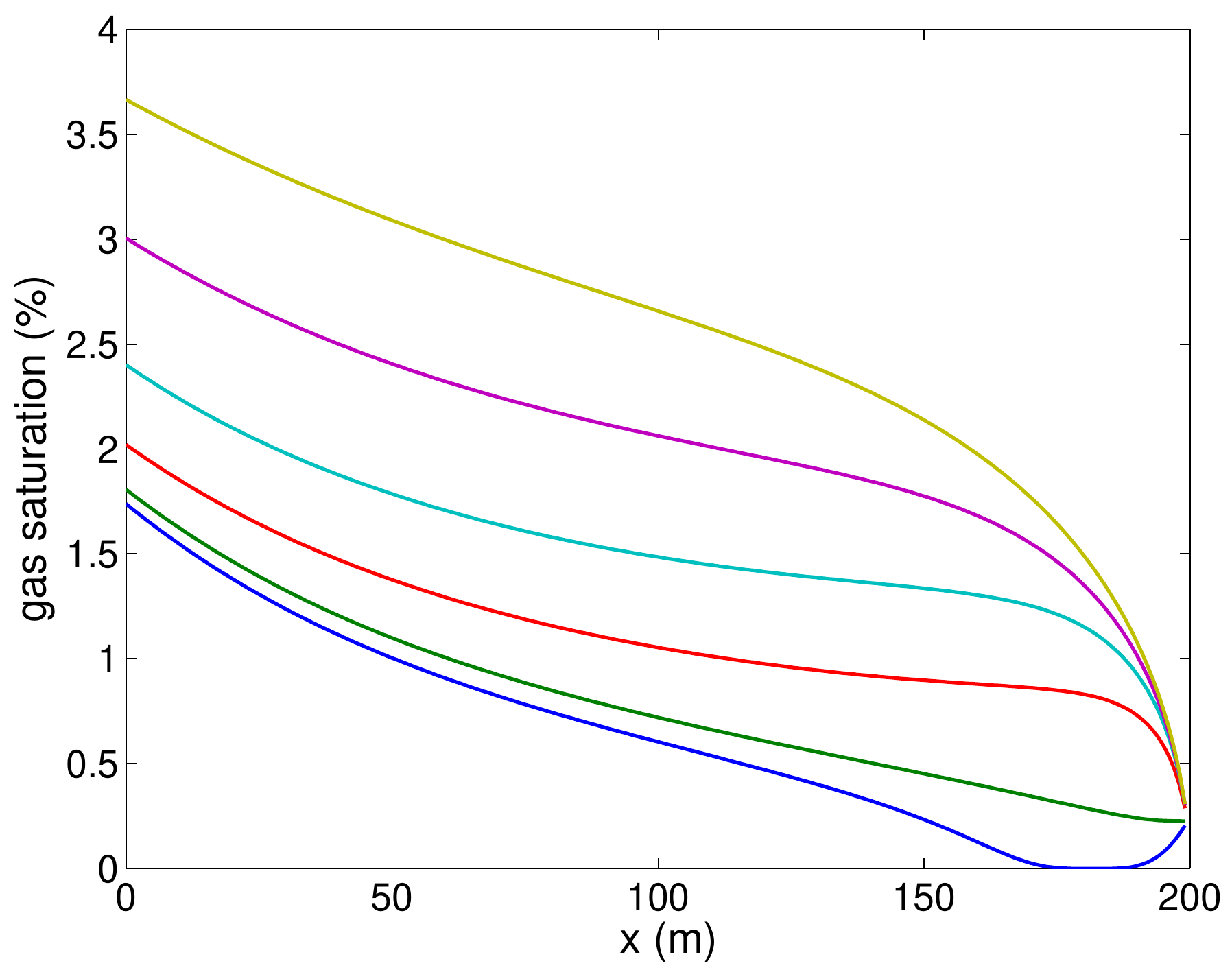}
    \includegraphics[width=.45\textwidth]{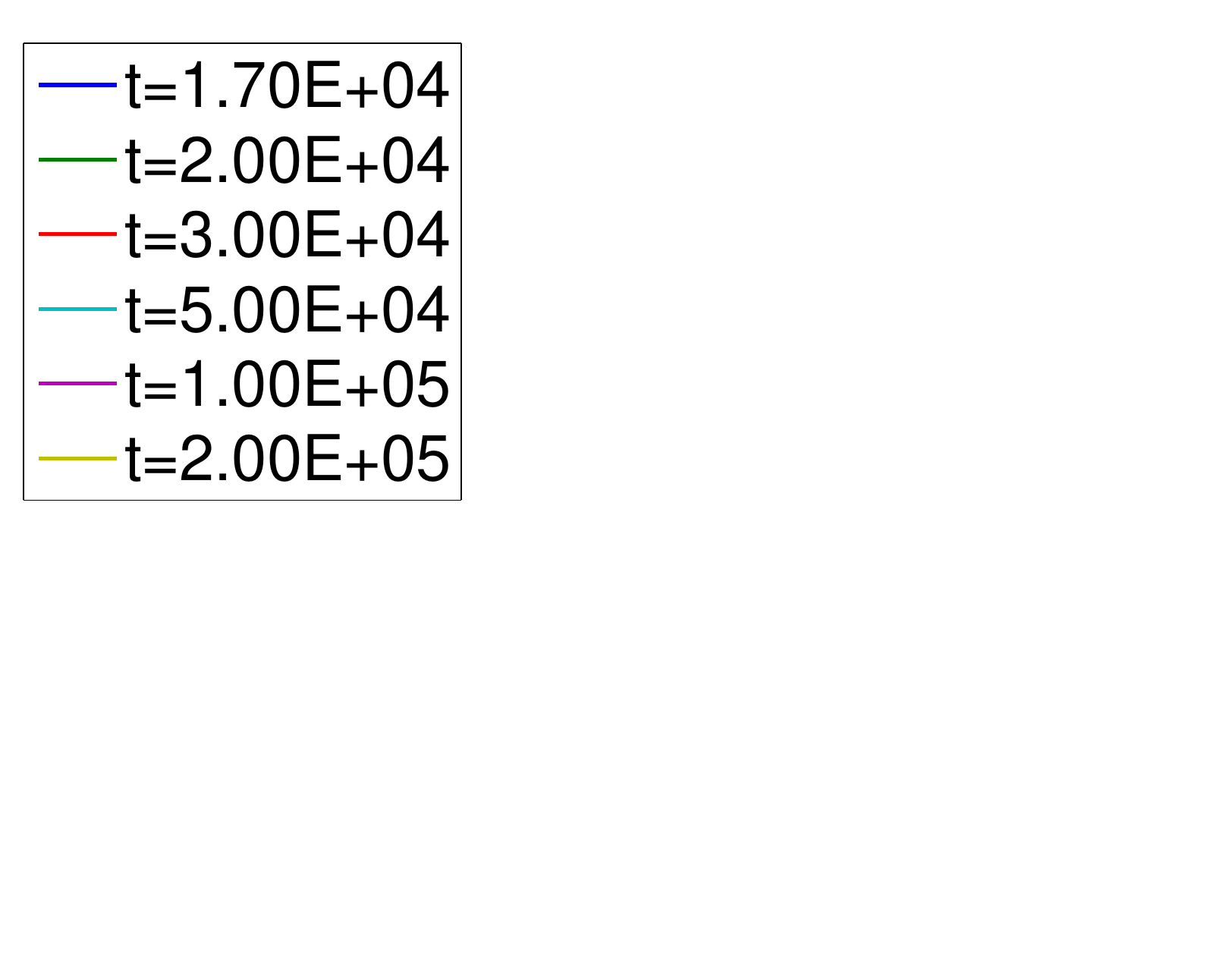}
    \caption{Test case number 2; $L_x= L_2=200$ m:  Time evolution of the dissolved hydrogen
    molar density ($\rho_l^h/M^h$) (top right),
   $p_l $(top left) and $S_g$ (bottom) profiles ;
    during the six last time steps.}
    \label{fig:ct2_2}
\end{figure}

 The main steps of the corresponding simulation are presented in Figures~\ref{fig:ct2_1}
 and~\ref{fig:ct2_2}where are presented
 the liquid
pressure $p_l$, the dissolved hydrogen molar density ( equal to
$\rho_l^h/M^h$) and the gas saturation $S_g$ profiles at different
times.

\begin{table}[htb] \centering \small
    \begin{tabular}{|c||c|rl||c|rl|}
        \hline
        Boundary conditions & \multicolumn{3}{c||}{Porous medium} & \multicolumn{3}{c|}{Others}\\
        \cline{2-7}
        Initial condition & Param. & \multicolumn{2}{c||}{Value} & Param. & \multicolumn{2}{c|}{Value}\\
        \hline
                $\phi^w\cdot\nu=0$ on $\Gamma_{imp}$ &
                $k$&$5\;10^{-20}$&\hspace{-1.5ex}$m^2$&
                $L_x$&$200$&\hspace{-1.5ex}$m$\\
                $\phi^h\cdot\nu=0$ on $\Gamma_{imp}$ &
                $\Phi$&$0.15$&$\hspace{-1.5ex}(-)$&
                $L_y$&$20$&\hspace{-1.5ex}$m$\\
                $\phi^w\cdot\nu=0$ on $\Gamma_{in}$ &
                $P_r$&$2\;10^6$&\hspace{-1.5ex}$Pa$&
                $L_1$&$0$&$m$\\
                $\phi^h\cdot\nu=\mathcal{Q}^h$  on $\Gamma_{in}$ &
                $n$&$1.49$&$\hspace{-1.5ex}(-)$&
                $p_{l,out}$&$10^6$&\hspace{-1.5ex}$Pa$\\
                $p_l=p_{l,out}$ on $\Gamma_{out}$ &
                $S_{l,res}$&$0.4$&$\hspace{-1.5ex}(-)$&
                $p_{g,out}$&$1.1\;10^6$&\hspace{-1.5ex}$Pa$\\
                $\rho_l^h=C_hp_{g,out}$ on $\Gamma_{out}$ &
                $S_{g,res}$&$0$&$\hspace{-1.5ex}(-)$&
                $\mathcal{Q}^h$&$55.7$&\hspace{-1.5ex}$mg/m^2/year$\\
                $p_l(t=0)=p_{l,out}$ in $\Omega$ &
                &&&
                &&\\
                $\rho_l^h(t=0)=C_hp_{g,out}$ in $\Omega$ &
                &&&
                &&\\
        \hline
    \end{tabular}
   \caption{Numerical Test case number 2: Boundary and Initial Conditions; porous medium characteristics and domain geometry.
   $\phi^w$ and $\phi^h$ are denoting respectively the water and hydrogen flux. }  \label{tab:ct2}
\end{table}

At the beginning, up to $t<1400$ years, the two phases are present
in the whole domain (see time $t=500$ years on
Figure~\ref{fig:ct2_1}). The permanent injection of hydrogen
increases both the two phase pressures and the gas saturation in the
vicinity of $\Gamma_{in}$. The local gas saturation drop is due to
the difference in mobilities between the two phases: the lower
liquid mobility leads to a bigger liquid pressure increase, compared
to the  gas pressure increase; which is finally  producing a
capillary pressure drop (according to definition (\ref{capillary}),
see Figure~\ref{P-C}), and creating  a water saturated zone. At time
$t=1400$ years, the gas phase starts to disappear in some region of
the porous domain (see time $t=1500$ years, in Figure
~\ref{fig:ct2_2}) .

Then, a saturated liquid region ($S_g=0$)will exist until time
$t=17\;000$ years (see  Figure ~\ref{fig:ct2_1}); and during this
period of time, the saturated region is pushed by the injected
Hydrogen, from $\Gamma_{in}$ to $\Gamma_{out}$.

After the time $t=17 \; 000$ years, due to the Dirichlet conditions
imposed on $\Gamma_{out}$, the liquid saturated  region disappears
and all together the phases pressure  and the gas saturation
 are growing in the whole domain (see the time $t=20 \;000$ years in
Figure~\ref{fig:ct2_2}).

Finally the liquid pressure reaches its maximum  at time $t=20
\;000$ years and then decreases in the whole domain (see the
Figure~\ref{fig:ct2_2}).  This is caused, like in the numerical test
case number 1, by  the evolution of the system
 towards a stationary state which is  characterized by a zero water component flow.


\subsection{Numerical Test number  3}
\label{Sec: 4.3}
 The geometry and the data of this numerical test
are given in Figure~\ref{fig:geom1D} and Table~\ref{tab:ct3} .
Like in the Numerical Test number 2, a constant
flux of hydrogen is imposed on the input boundary,
 $\Gamma_{in}$, while Dirichlet conditions $p_l=p_{l,out}$, $\rho_l^h=0$ are given
 on $\Gamma_{out}$, in order to have only the liquid phase on this part of the boundary. The initial conditions, $p_l=p_{l,out}$ and
 $\rho_l^h=0$, are uniform on all the domain, and correspond to a porous domain initially saturated with pure
 water.
Contrary to the two first numerical tests, the porous domain is non homogeneous, there are two different porous subdomains
 $\Omega_1$ and $\Omega_2$; $L_x= 200$ m, $L_1= 20$ m and $L_2= 180$ m.\\

\begin{table}[htb] \centering \small
    \begin{tabular}{|c||c|rl||c|rl|}
        \hline
        Boundary conditions & \multicolumn{3}{c||}{Porous medium} & \multicolumn{3}{c|}{Other}\\
        \cline{2-7}
        initial condition & Param. & \multicolumn{2}{c||}{Value on $\Omega_1$} & Param. & \multicolumn{2}{c|}{Value}\\
        \hline
        $\phi^w\cdot\nu=0$ on $\Gamma_{imp}$ &
        $k$&$10^{-18}$&\hspace{-1.5ex}$m^2$&
        $L_x$&$200$&\hspace{-1.5ex}$m$\\
        $\phi^h\cdot\nu=0$ on $\Gamma_{imp}$ &
        $\Phi$&$0.3$&\hspace{-1.5ex}$(-)$&
        $L_y$&$20$&\hspace{-1.5ex}$m$\\
        $\phi^w\cdot\nu=0$ on $\Gamma_{in}$ &
        $P_r$&$2\;10^6$&\hspace{-1.5ex}$Pa$&
        $L_1$&$20$&\hspace{-1.5ex}$m$\\
        $\phi^h\cdot\nu=\mathcal{Q}^h$ on $\Gamma_{in}$ &
        $n$&$1.54$&\hspace{-1.5ex}$(-)$&
        $p_{l,out}$&$10^6$&\hspace{-1.5ex}$Pa$\\
        $p_l=p_{l,out}$ on $\Gamma_{out}$ &
        $S_{l,res}$&$0.01$&\hspace{-1.5ex}$(-)$&
        $\mathcal{Q}^h$&$5.57$&\hspace{-1.5ex}$mg/m^2/year$\\
        $\rho_l^h=0$ on $\Gamma_{out}$ &
        $S_{g,res}$&$0$&\hspace{-1.5ex}$(-)$&
        $ $&$ $&$ $\\
        \cline{2-4}
        $p_l(t=0)=p_{l,out}$ on $\Omega$ &
        Param.&\multicolumn{2}{c||}{Value on $\Omega_2$}&
        &&\\
        \cline{2-4}
        $\rho_l^h(t=0)=0$ on $\Omega$ &
        $k$&$5\;10^{-20}$&\hspace{-1.5ex}$m^2$&
        &&\\
        &
        $\Phi$&$0.15$&\hspace{-1.5ex}$(-)$&
        &&\\
        &
        $P_r$&$15\;10^6$&\hspace{-1.5ex}$Pa$&
        &&\\
        &
        $n$&$1.49$&\hspace{-1.5ex}$(-)$&
        &&\\
        &
        $S_{l,res}$&$0.4$&\hspace{-1.5ex}$(-)$&
        &&\\
        &
        $S_{g,res}$&$0$&\hspace{-1.5ex}$(-)$&
        &&\\        %
        \hline
    \end{tabular}
    \caption{Numerical Test case number 3: Boundary and Initial Conditions; porous medium characteristics and domain geometry.
    $\phi^w$ and $\phi^h$ are denoting respectively the water and hydrogen flux.} \label{tab:ct3}
\end{table}

\begin{figure}[hbtp]
    \includegraphics[width=.45\textwidth]{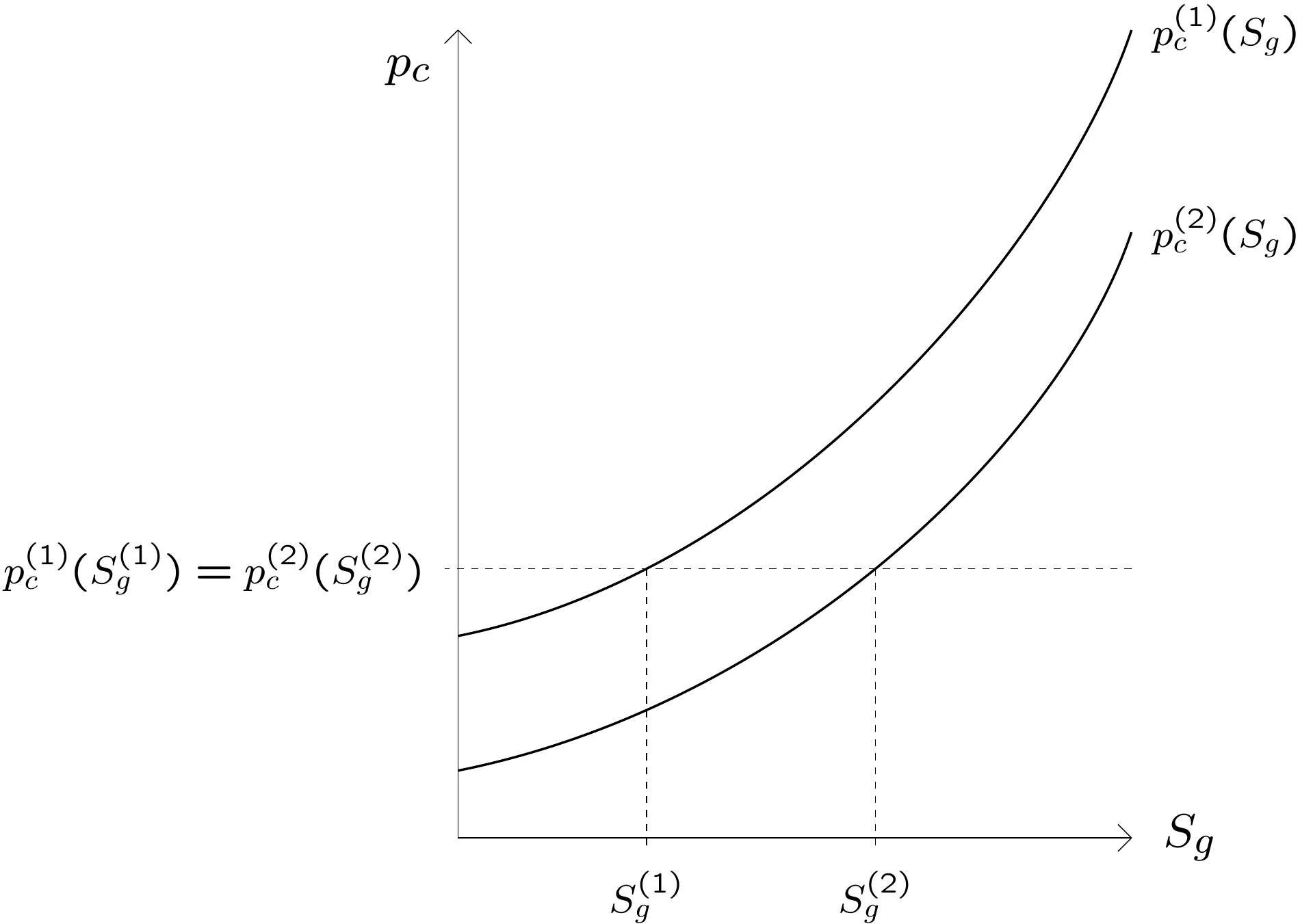}
    \caption{Saturation discontinuity at the
    interface of two materials with different capillary pressure curves; test case number 3.}
    \label{fig:Pc_disc_t}
\end{figure}

\begin{figure}[hbtp]
    \includegraphics[width=.45\textwidth]{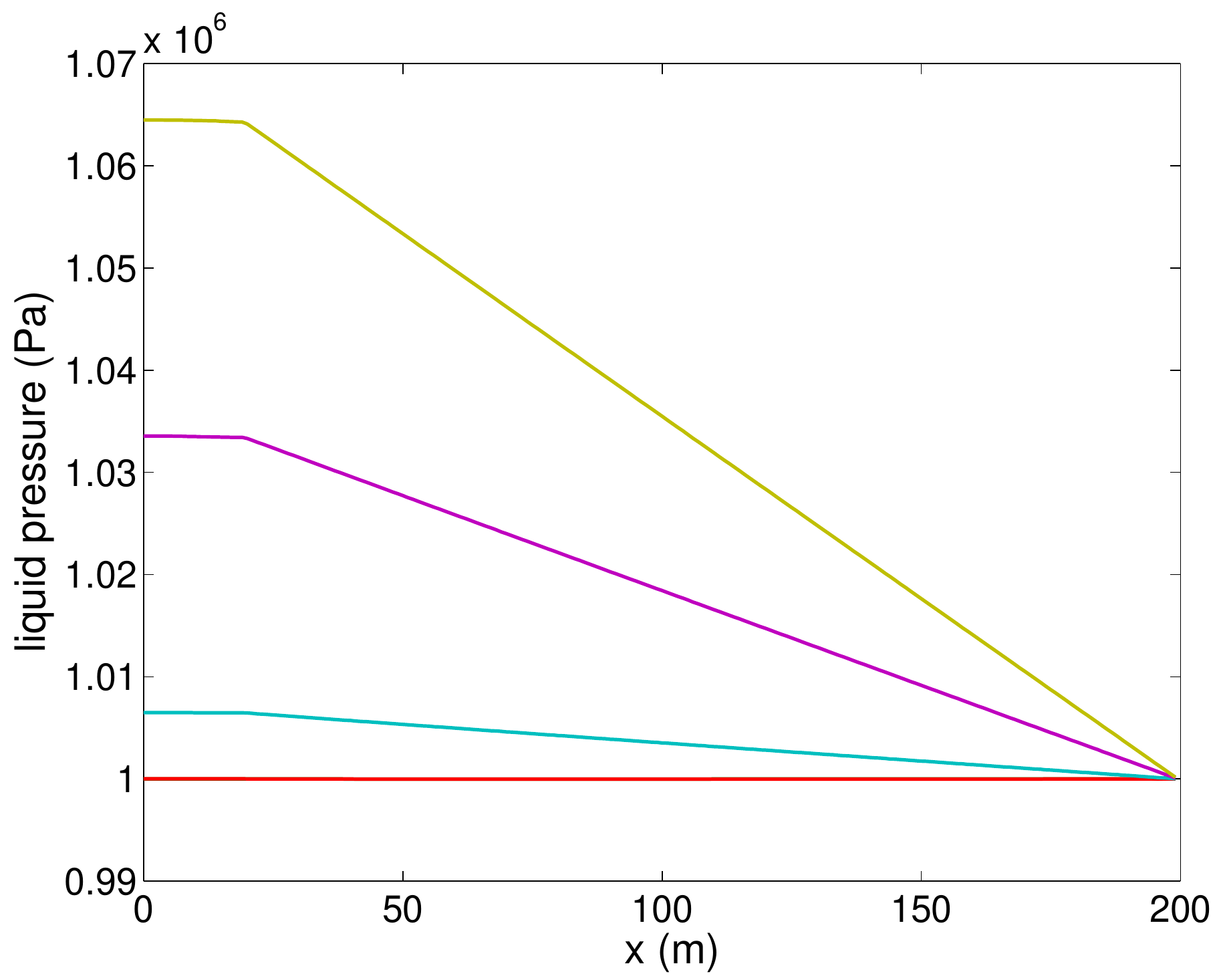}
    \includegraphics[width=.45\textwidth]{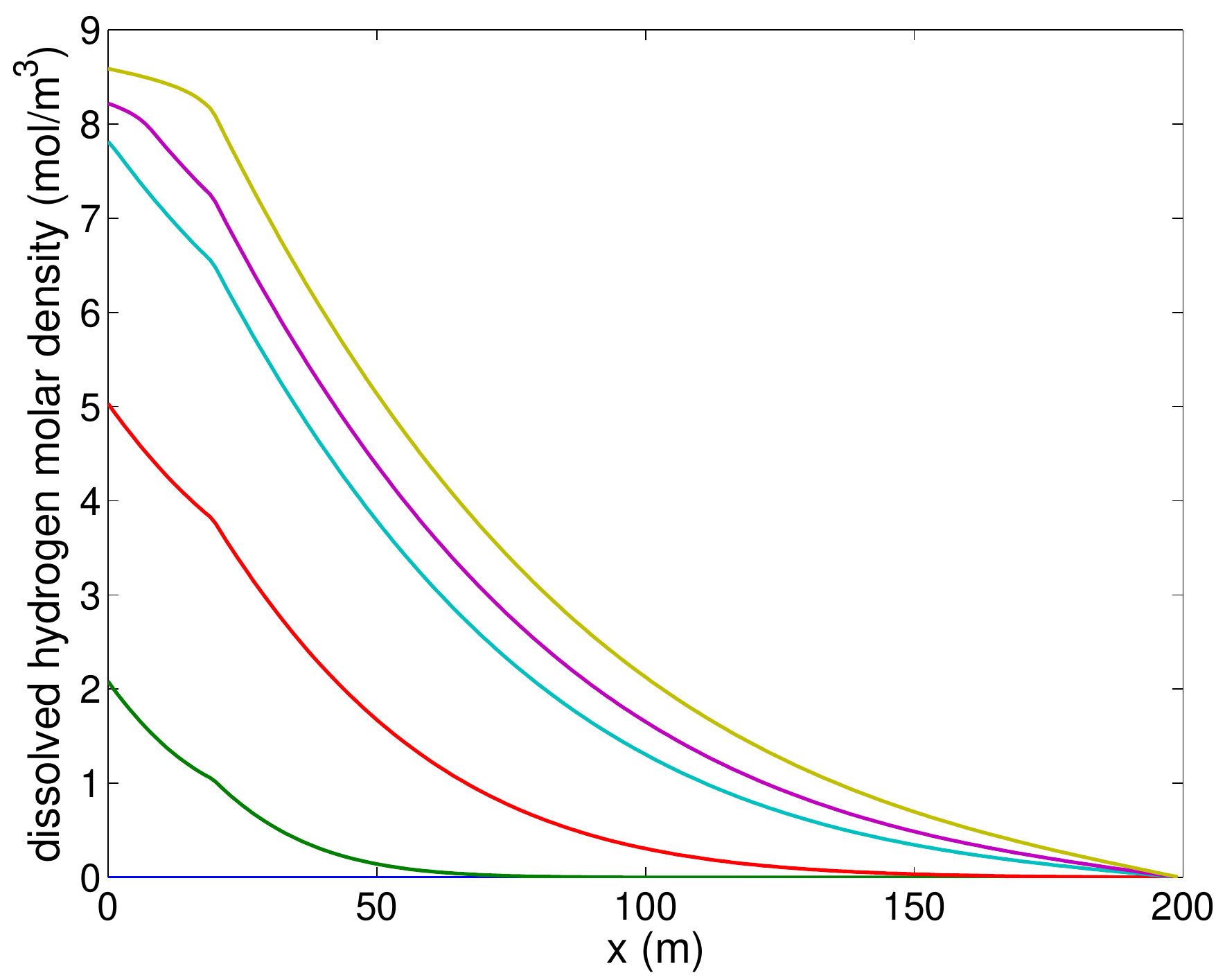}
    \includegraphics[width=.45\textwidth]{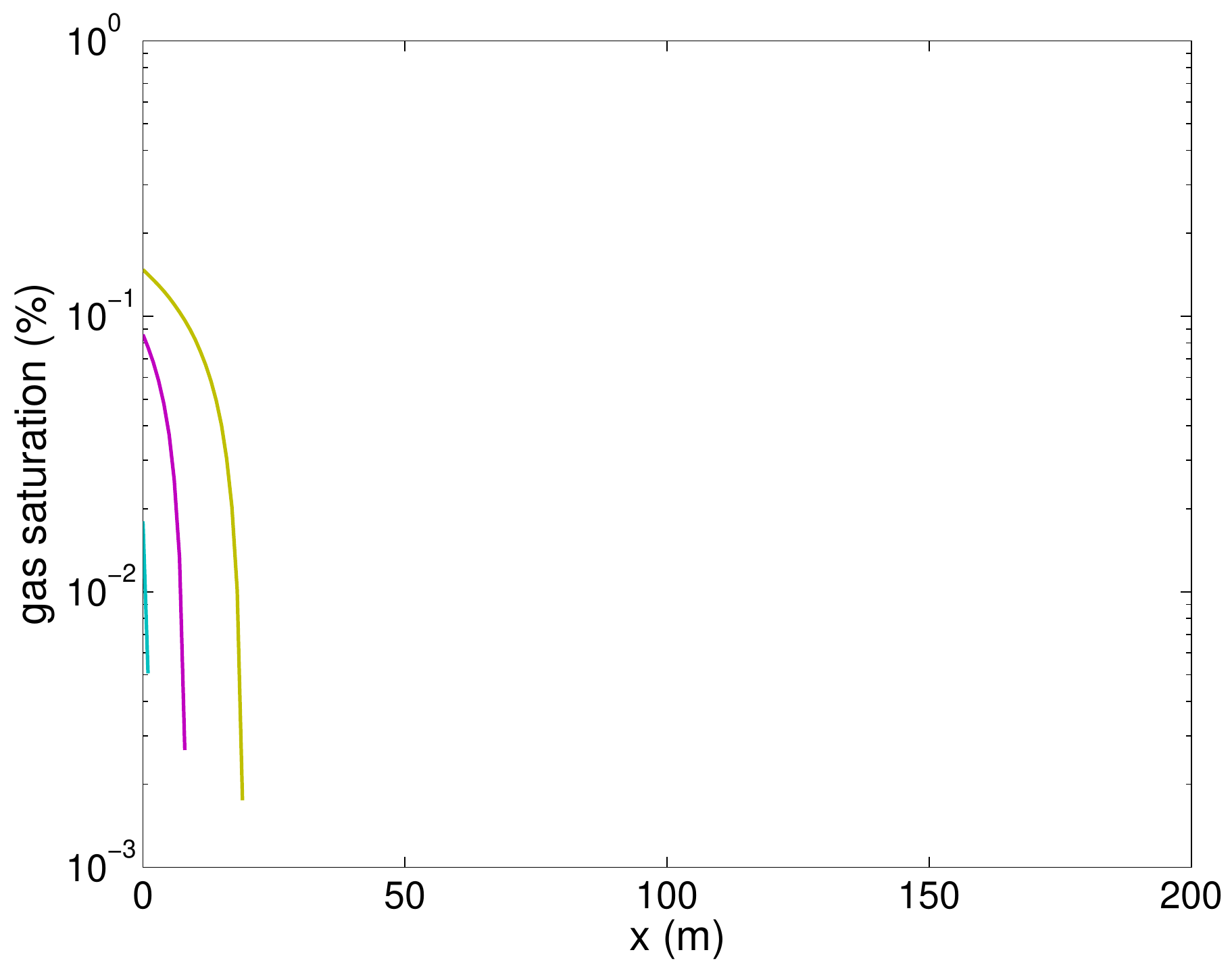}
    \includegraphics[width=.45\textwidth]{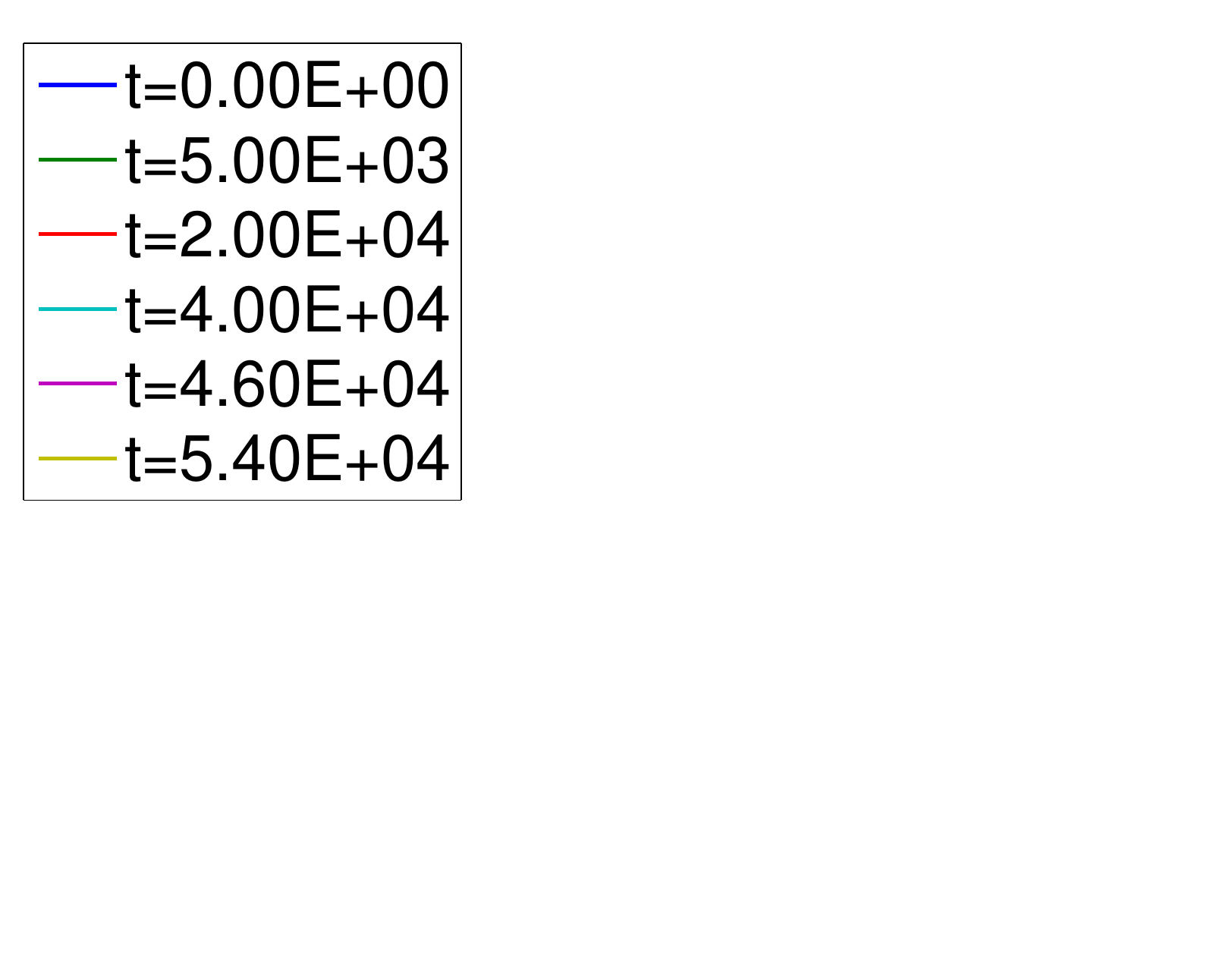}
    \caption{Test case number 3; $L_x= 200$ m, $L_1= 20$m:  Time evolution of the dissolved hydrogen
    molar density ($\rho_l^h/M^h$) (top right), $p_l$ (top left) and $S_g$ (bottom) profiles ;
    during the first time steps. All the $S_g$ curves go to zero( although this cannot be seen using a logarithmic scale) .}
    \label{fig:ct3_1}
\end{figure}

\begin{figure}[hbtp]
    \includegraphics[width=.45\textwidth]{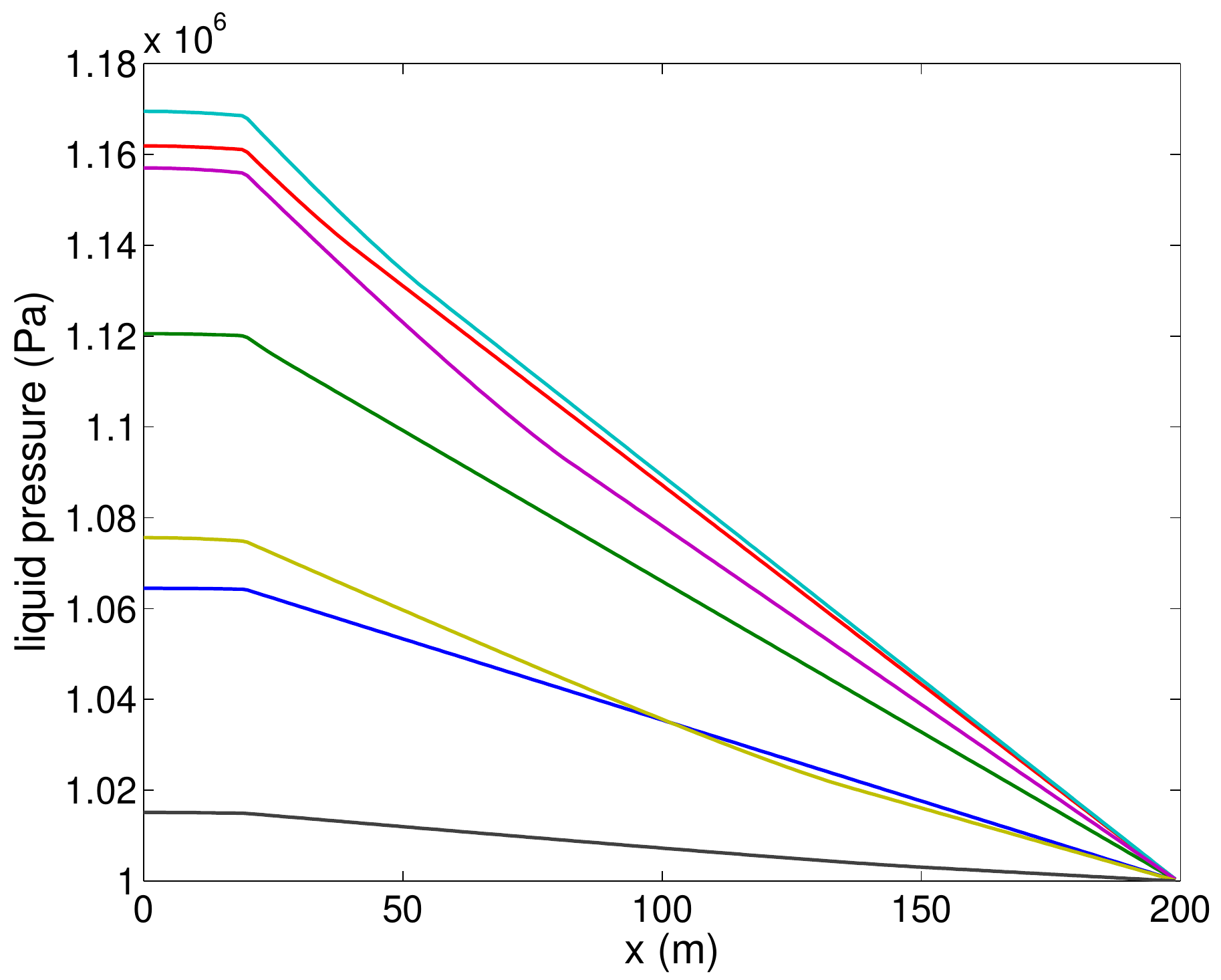}
    \includegraphics[width=.45\textwidth]{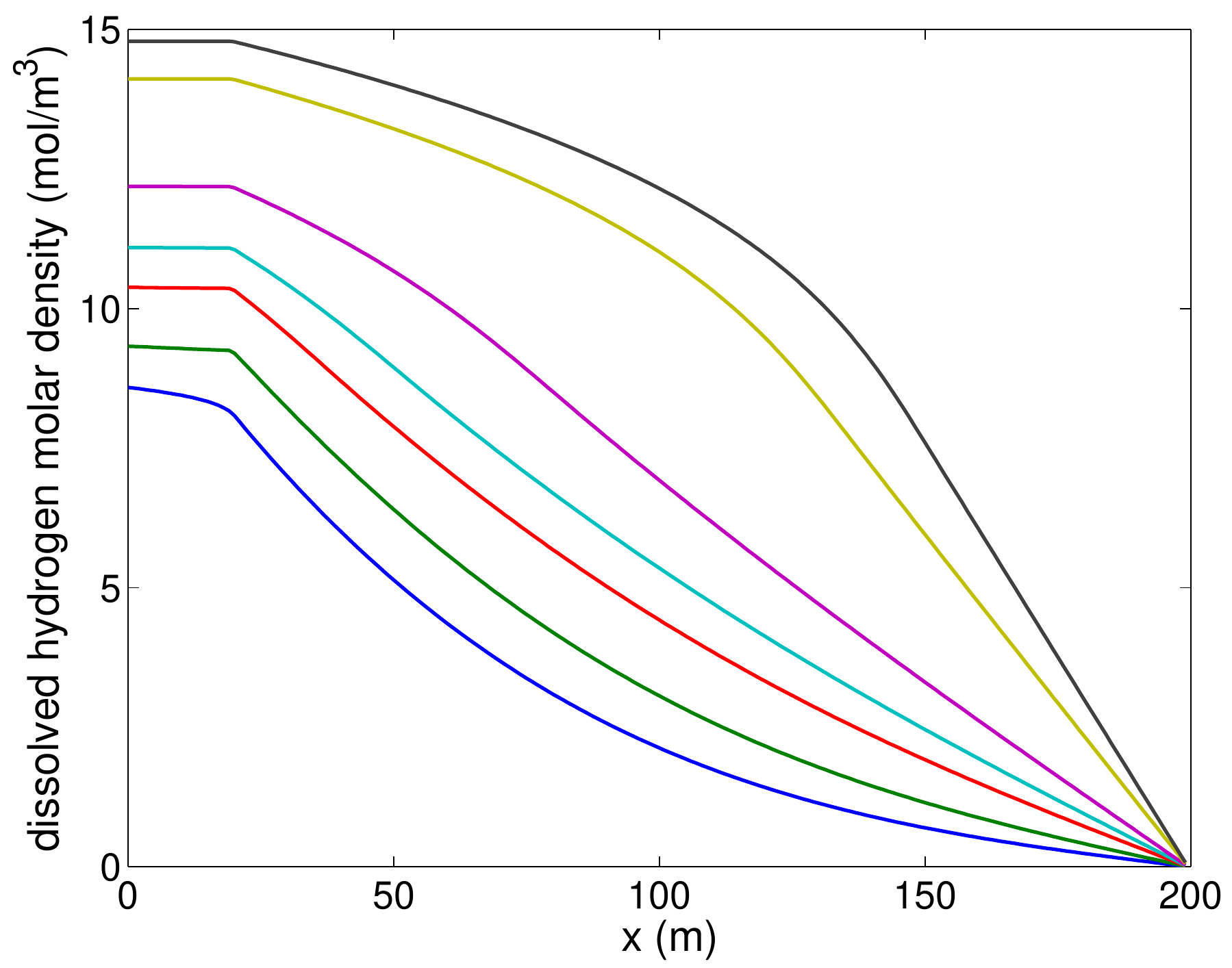}
    \includegraphics[width=.45\textwidth]{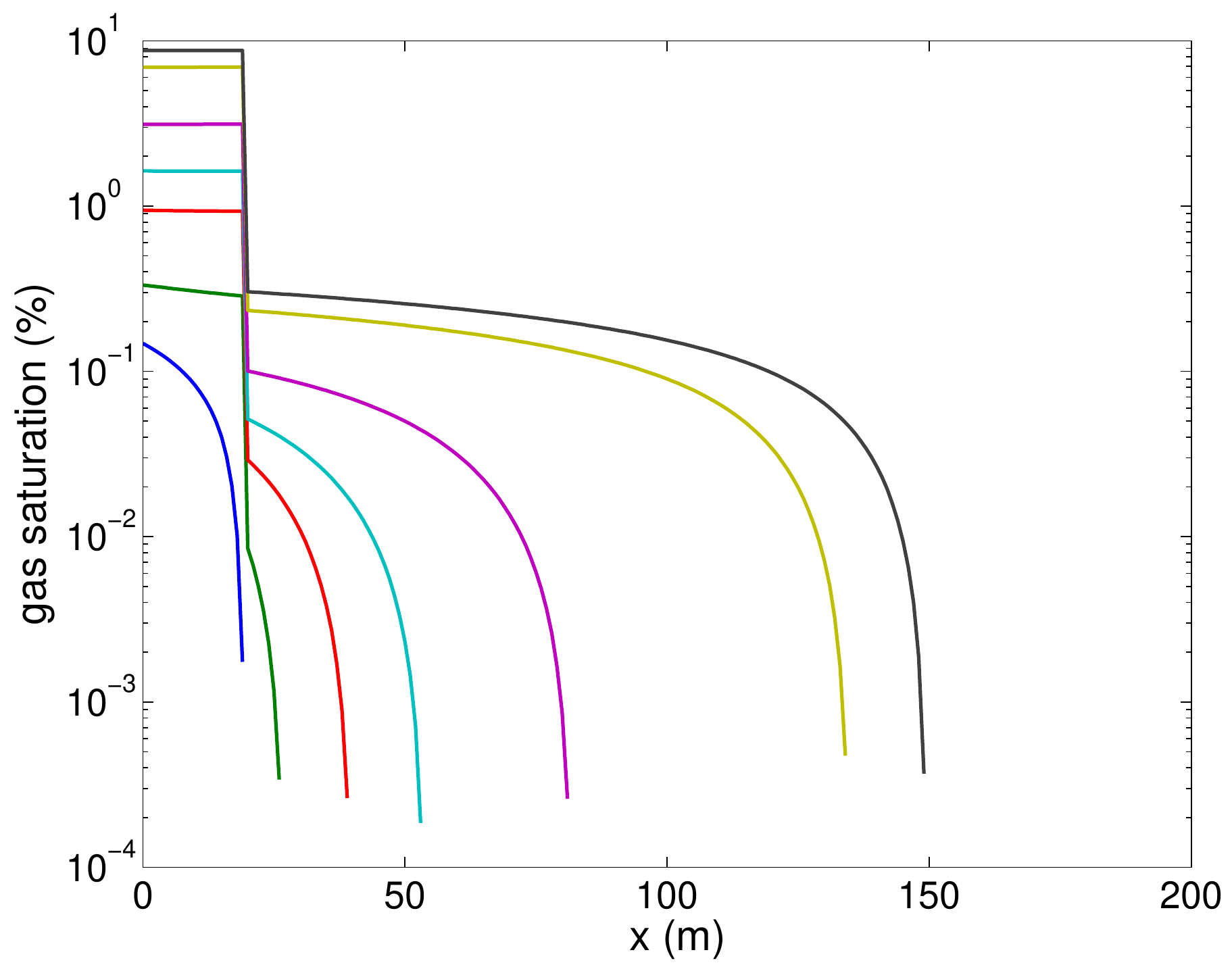}
    \includegraphics[width=.45\textwidth]{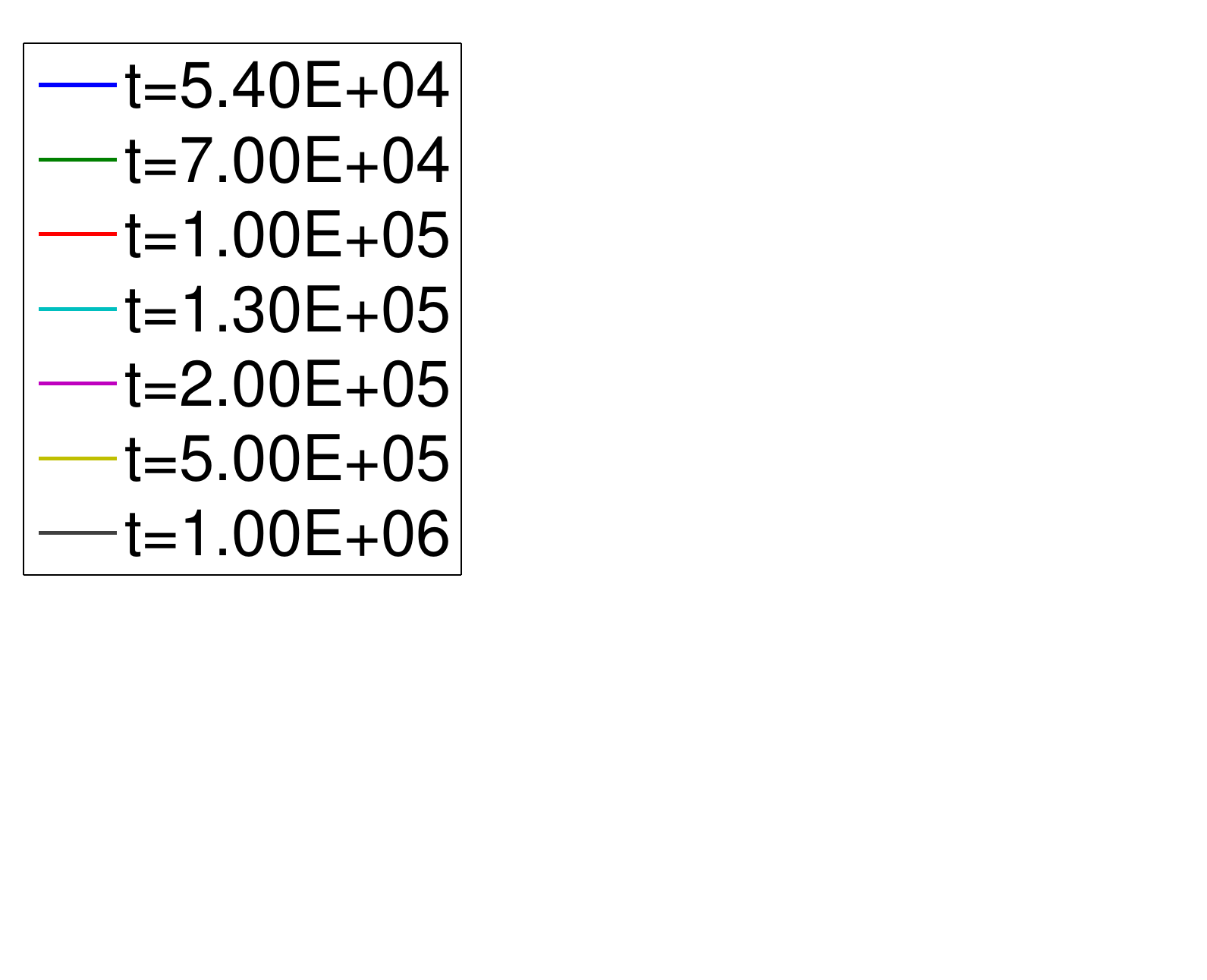}
    \caption{Test case number 3; $L_x= 200$ m, $L_1= 20$m:  Time evolution of the dissolved hydrogen
    molar density ($\rho_l^h/M^h$)(top right), $p_l$ (top left) and $S_g$ (bottom) profiles ;
    during the last seven time steps. All the $S_g$ curves go to zero( although this cannot be seen using a logarithmic scale) .}
    \label{fig:ct3_2}
\end{figure}

The simulation time of this test case is  $T=10^6$ years;the
discretization space mesh is 1 m; the time step is  $10^2$ years at
the beginning and grows up to
 $2\cdot 10^4$ years in the end of the simulation (see Table\ref{tab:NUM}).\\
Figures~\ref{fig:ct3_1} and~\ref{fig:ct3_2} represent the liquid
pressure $p_l$, the dissolved hydrogen molar density ( equal to
$\rho_l^h/M^h$) and the gas saturation $S_g$ profiles at different
times.

 The main difference from the previous  simulations (which
were in a homogeneous porous domain)  is the gas saturation
discontinuity, staying on the porous domain interface $x=20 \; $m;
and due to the height of this saturation jump ,
  we had to use a logarithm scale for presenting the gas saturation $S_g$ profiles .

There are four main steps :
\begin{itemize}
    \item From 0  to $3.8\cdot 10^4$ years both the gas saturation and the liquid pressure stay constant in the whole domain
        while the hydrogen injection on the left side $\Gamma_{in}$ of the domain  increases the hydrogen density level .
    \item From $3.8\cdot 10^4$ to $5.4\cdot 10^4$ years  both the liquid pressure  and  the hydrogen density
    are increasing
        in the whole domain. The gas  start  to expanding  from the left side  of the domain $\Gamma_{in}$. The
        saturation front is moved towards the porous media discontinuity, at $x=20\;m$, which is  reached at
        $t=5.4\cdot 10^4$ years; see Figures~\ref{fig:ct3_1}.
    \item From  $5.4\cdot10^4$ years to $1.3\cdot 10^5$ years, see Figures~\ref{fig:ct3_2}, the saturation front
       has crossed the medium discontinuity at  $x=20\;m$  and, from now, all the saturation profiles will have
        a discontinuity at $x=20\;m$.
    \item From $1.3\cdot 10^5$ years to $10^6$ years, see Figures~\ref{fig:ct3_2}, both the  hydrogen
     density and the gas  saturation  keep
        growing while the liquid pressure decreases towards zero on the entire domain. The gas saturation front
        keeps  moving to the right, pushed by the injected gas, up
        to  $x\approx 150\;m$ at $10^6$ years.
\end{itemize}

Until the saturation front reaches the interface between the two
porous media, for ($t=5.4\cdot 10^4$ years),
 appearance and evolution of both the gas phase and the unsaturated zone are identical to
what was happening  in the test case 1 (with a homogeneous porous domain) during the period of gas injection:
the dissolved hydrogen is accumulating at the entrance until  the liquid phase becomes saturated,at time
 ($t>3.8\:10^4$ years), letting the gas phase to appear.

When the saturation front crosses the interface between the two
porous subdomains (at $x=20\;m$ and $t=5.4\cdot 10^4$ years), the
gas saturation  is strictly positive on both sides of this interface
and  the caplllary pressure curves being different on each side( see
Table~\ref{tab:ct3})  forces  the saturation to be discontinuous for
preserving
 the capillary pressure continuity on the interface.
 The capillary pressure
continuity at the interface imposes to $p_c^{(1)}$, the Capillary
Pressure in $\Omega_1$, and to $p_c^{(2)}$, the Capillary Pressure
in $\Omega_2$, to be equal on this interface.
 $p_c^{(1)} =p_c^{(2)} $ is satisfied only if there  are two  different  saturations,on each interface side $S_g^{(1)}$,
 and $S_g^{(2)}$:
  $p_c^{(1)}(S_g^{(1)})=p_c^{(2)}(S_g^{(2)})$ ; see Figure \ref{fig:Pc_disc_t}.

In the same way as in the numerical test case number 1, the system
tends to a stationary state


\subsection{Numerical Test number 4}
\label{Sec: 4.4}
 This last numerical test  is different from  all
the precedent ones; it intends to be a simplified representation of
  what happens when an unsaturated porous block is placed  within a water saturated porous structure.
The challenge is then: how   the mechanical balance will be restored
in a homogeneous porous domain , which was initially out of
equilibrium, i.e.  with a jump in the initial phase pressures?\\

 The initial liquid pressure  is the same in the entire  porous
 domain;
$\Omega$ ,$p_{l,1}=p_{l,2}$ , and in the subdomain $\Omega_1$ the
initial condition, ($p_{l,1}=p_{g,1}$ in Table~\ref{tab:ct4},
corresponds to a liquid fully saturated state with a hydrogen
concentration reaching the gas appearance concentration threshold
($p_g=p_l$ and $\rho_{l}^{h}=C_{h}p_{g}$, see Figure~\ref{Rem:1}).
In the subdomain $\Omega_2$ the initial condition($p_{l,2}\neq
p_{g,2}$and $p_{g,2} \neq p_{g,1}$ )corresponds to a non saturated
state (see Table~\ref{tab:ct4}).\\
The porous block initial state is said out of equilibrium,
because:\\
if this initial state  was in equilibrium , in the two subdomains
$\Omega_1\;$ and $\Omega_2$, the local mechanical balance would have
made the pressures, of both the liquid and the gas phase, continuous
in the entire domain $\Omega\;$.

  For simplicity, we assume the porous medium domain $\Omega$ is
homogeneous and all the porous medium characteristics are the same
in the two subdomains $\Omega_1$ and $\Omega_2$, and corresponding
to concrete.

The system is then expected to evolve from this initial out of
equilibrium state towards a stationary state.

We should notice that, in order to see appearing  the  final
stationary state, in a reasonable period of time,
 we have shortened  the domain $\Omega$ ( $L_x = 1$m),
 taken the  porous media characteristics , and
set the final time of this simulation $T_{fin}$ at
$T_{fin}=10^6\;s\approx11.6\;$) days.  The complete set of data
 of this test case is given in Table~\ref{tab:ct4}.

\begin{table}[htb] \centering \small
    \begin{tabular}{|c||c|rl||c|rl|}
        \hline
        Boundary conditions & \multicolumn{3}{c||}{Porous medium} & \multicolumn{3}{c|}{Other}\\
        \cline{2-7}
        initial condition & Param. & \multicolumn{2}{c||}{Value} & Param. & \multicolumn{2}{c|}{Value}\\
        \hline
        $\phi^w\cdot\nu=0$ on $\partial\Omega$ &
        $k$&$10^{-18}$&\hspace{-1.5ex}$m^2$&
        $L_x$&$1$&\hspace{-1.5ex}$m$\\
        $\phi^h\cdot\nu=0$ on $\partial\Omega$ &
        $\Phi$&$0.3$&$\hspace{-1.5ex}(-)$&
        $L_y$&$0.1$&\hspace{-1.5ex}$m$\\
        $p_l(t=0)=p_{l,1}$ on $\Omega_1$ &
        $P_r$&$2\;10^6$&\hspace{-1.5ex}$Pa$&
        $L_1$&$0.5$&$m$\\
        $\rho_l^h(t=0)=C_hp_{g,1}$ on $\Omega_1$ &
        $n$&$1.54$&$\hspace{-1.5ex}(-)$&
        $p_{l,1}$&$10^6$&\hspace{-1.5ex}$Pa$\\
        $p_l(t=0)=p_{l,2}$ on $\Omega_2$ &
        $S_{l,res}$&$0.01$&$\hspace{-1.5ex}(-)$&
        $p_{g,1}$&$10^6$&\hspace{-1.5ex}$Pa$\\
        $\rho_l^h(t=0)=C_hp_{g,2}$ on $\Omega_2$ &
        $S_{g,res}$&$0$&$\hspace{-1.5ex}(-)$&
        $p_{l,2}$&$10^6$&\hspace{-1.5ex}$Pa$\\
        &
        &&&
        $p_{g,2}$&$2.5\;10^6$&\hspace{-1.5ex}$Pa$\\
        \hline
    \end{tabular}

    \caption{Data of the numerical test number 4 : boundary and initial conditions;domain geometry. The porous medium domain
    $\Omega$ is homogeneous, all the porous medium parameters are the same in the two subdomains $\Omega_1$ and $\Omega_2$;
     $\phi^w$ and $\phi^h$ are denoting respectively the water and hydrogen flux .}
     \label{tab:ct4}
\end{table}

\begin{figure}[hbtp]
    \includegraphics[width=.45\textwidth]{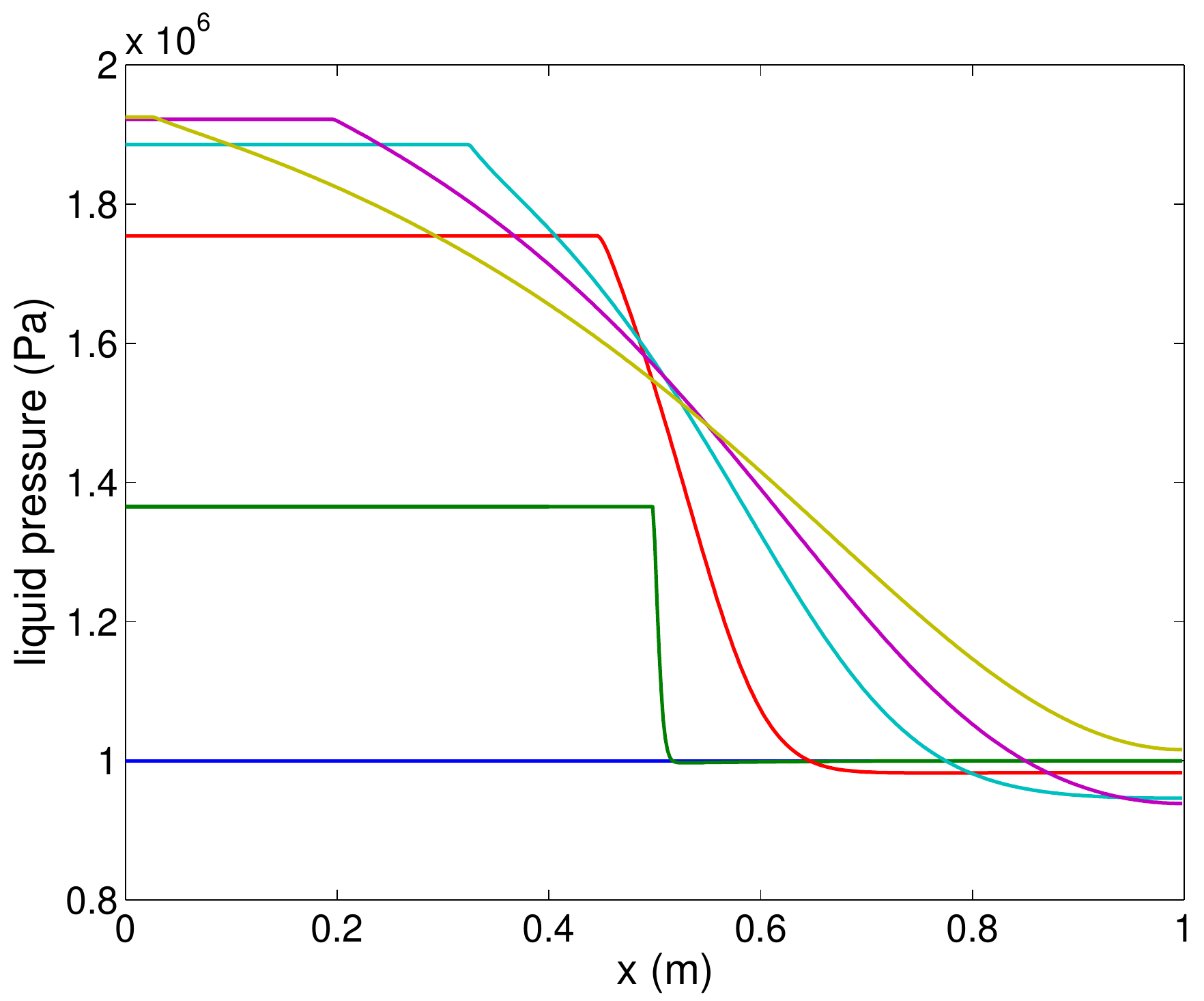}
    \includegraphics[width=.45\textwidth]{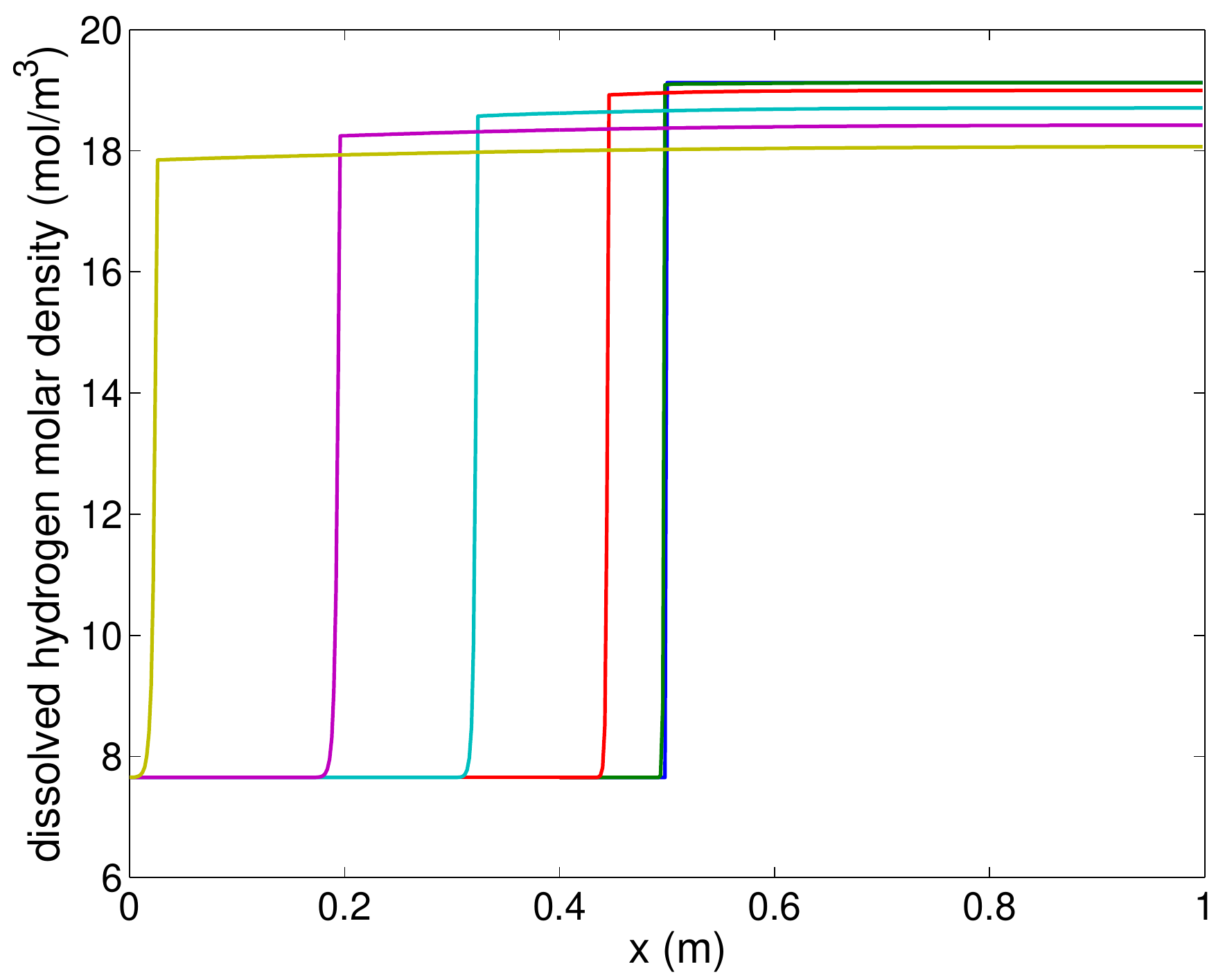}
    \includegraphics[width=.45\textwidth]{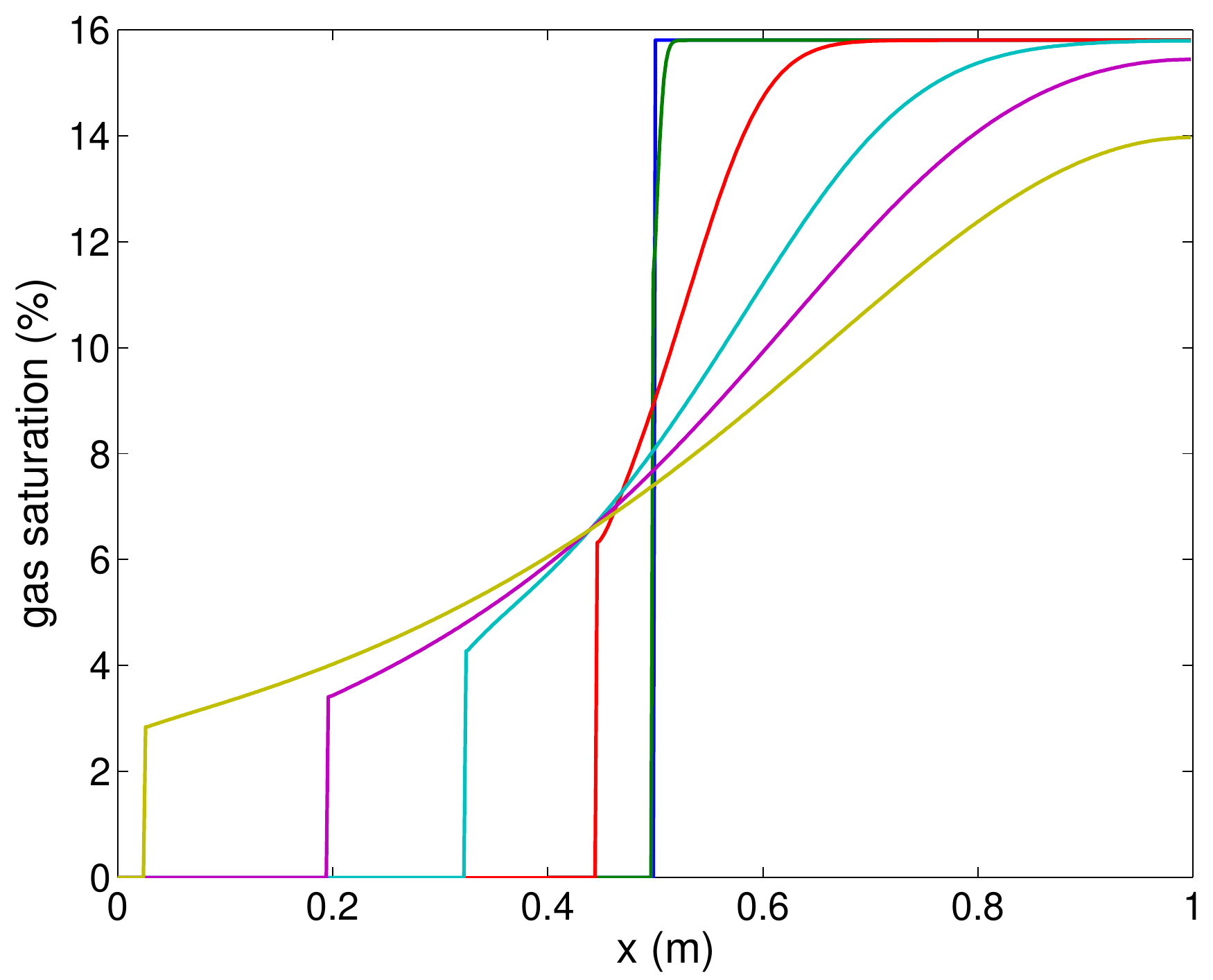}
    \includegraphics[width=.45\textwidth]{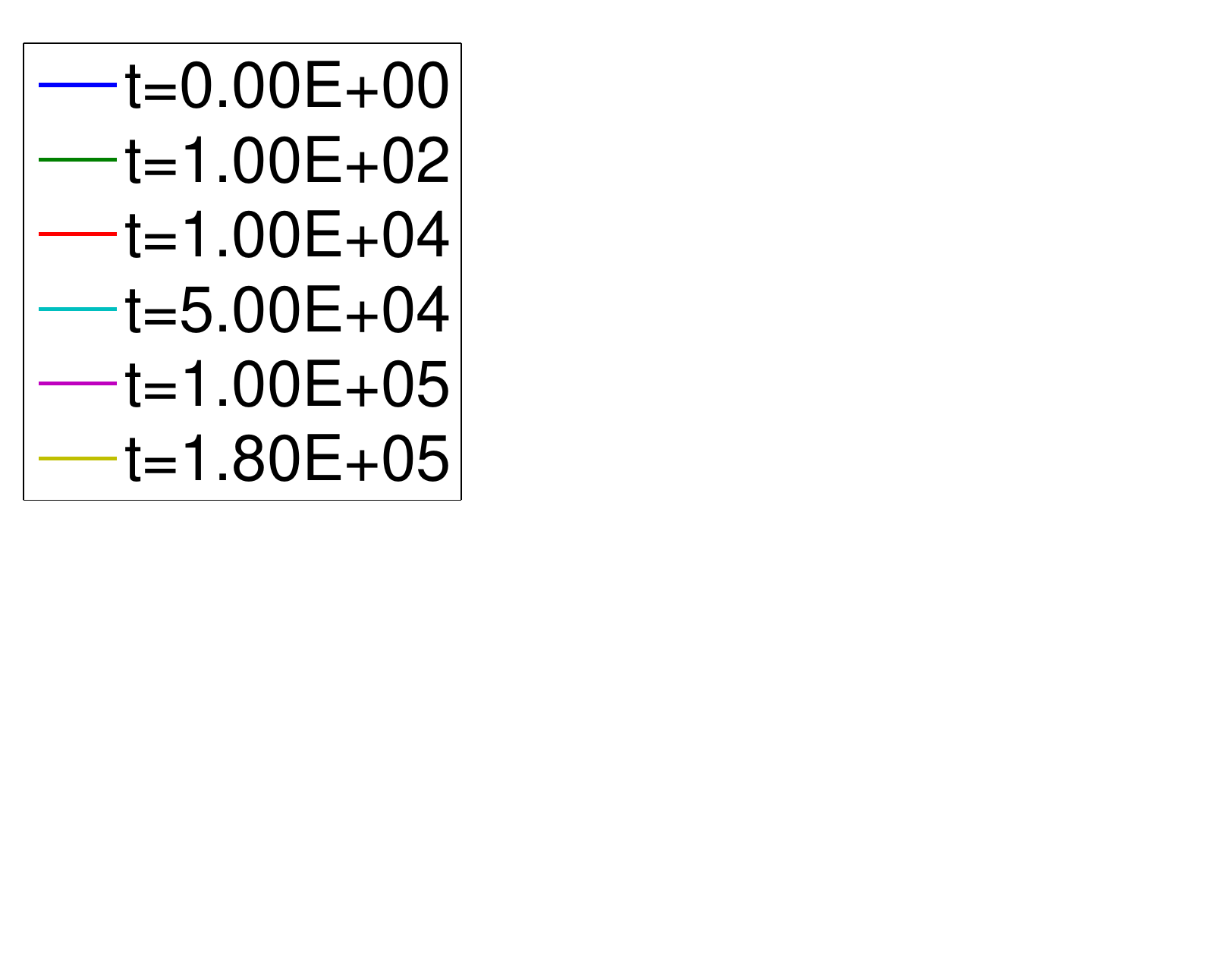}
    \caption{Numerical test case number 4, $L_x=1$ m, $L_1=0.5$ m: Time evolution of the dissolved hydrogen
    molar density ($\rho_l^h/M^h$) (top right), $p_l$ (top left)
    and $S_g$ (bottom) profiles; during the six first time steps.}
    \label{fig:ct4_1}
\end{figure}
\begin{figure}[hbtp]
    \includegraphics[width=.45\textwidth]{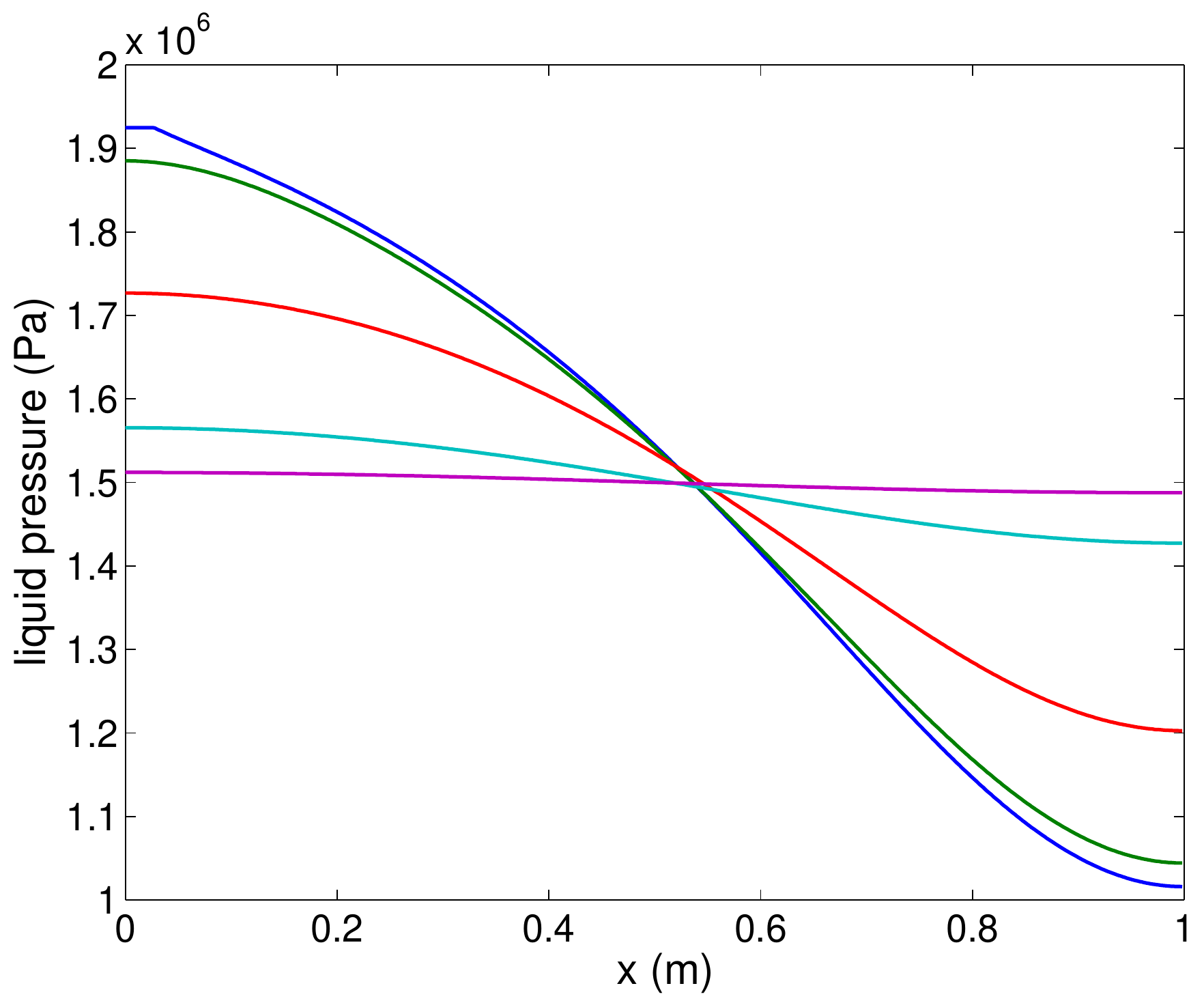}
    \includegraphics[width=.45\textwidth]{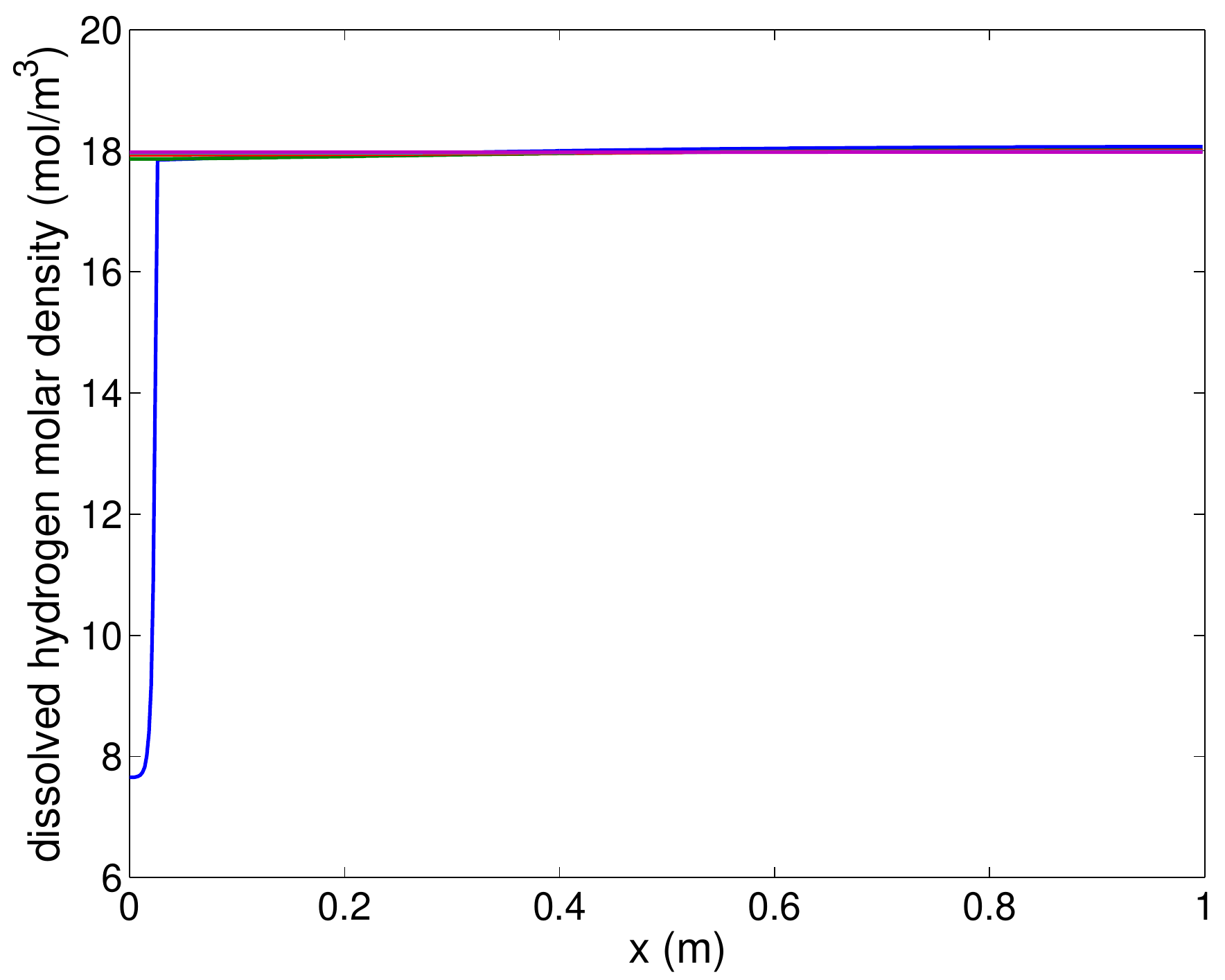}
    \includegraphics[width=.45\textwidth]{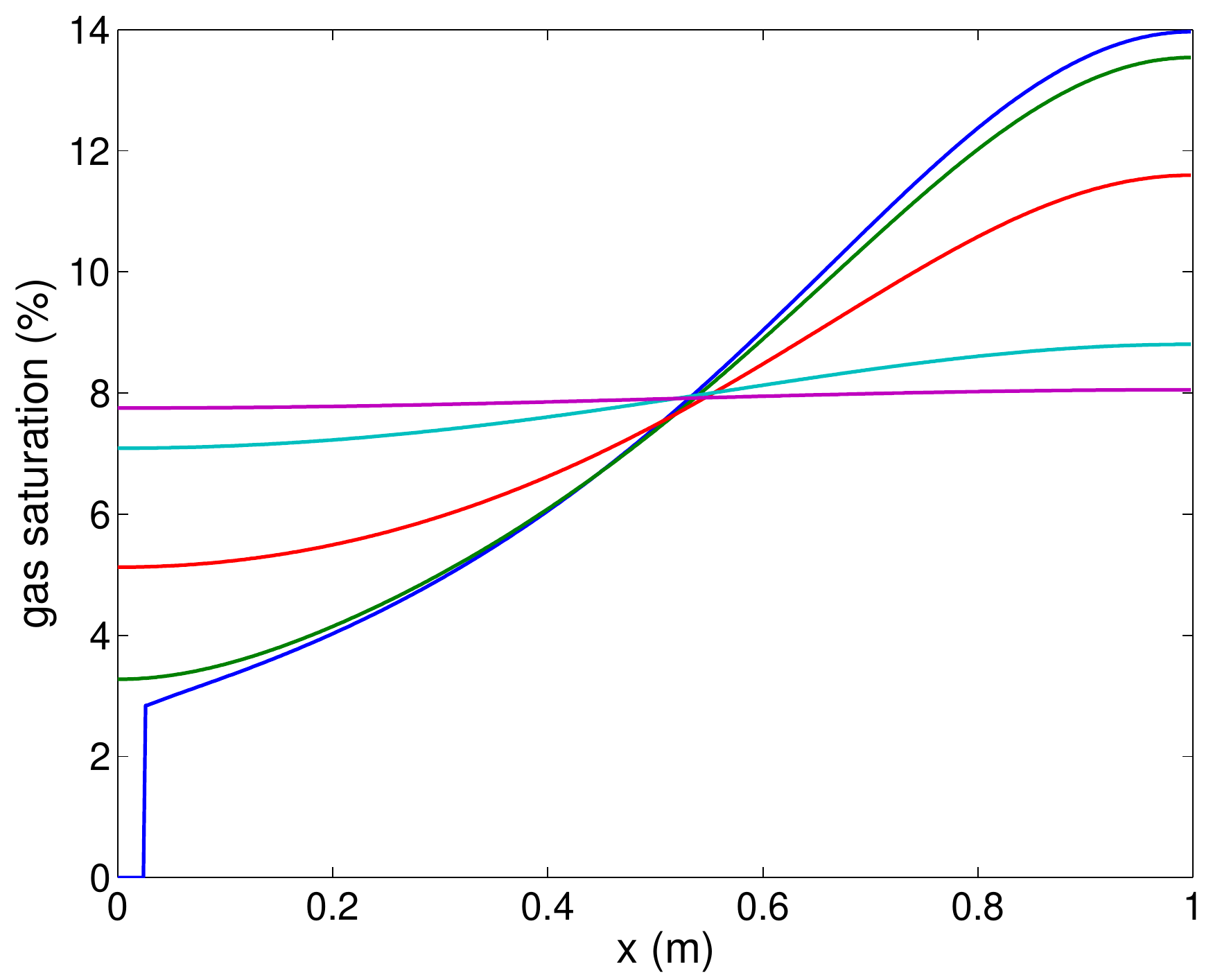}
    \includegraphics[width=.45\textwidth]{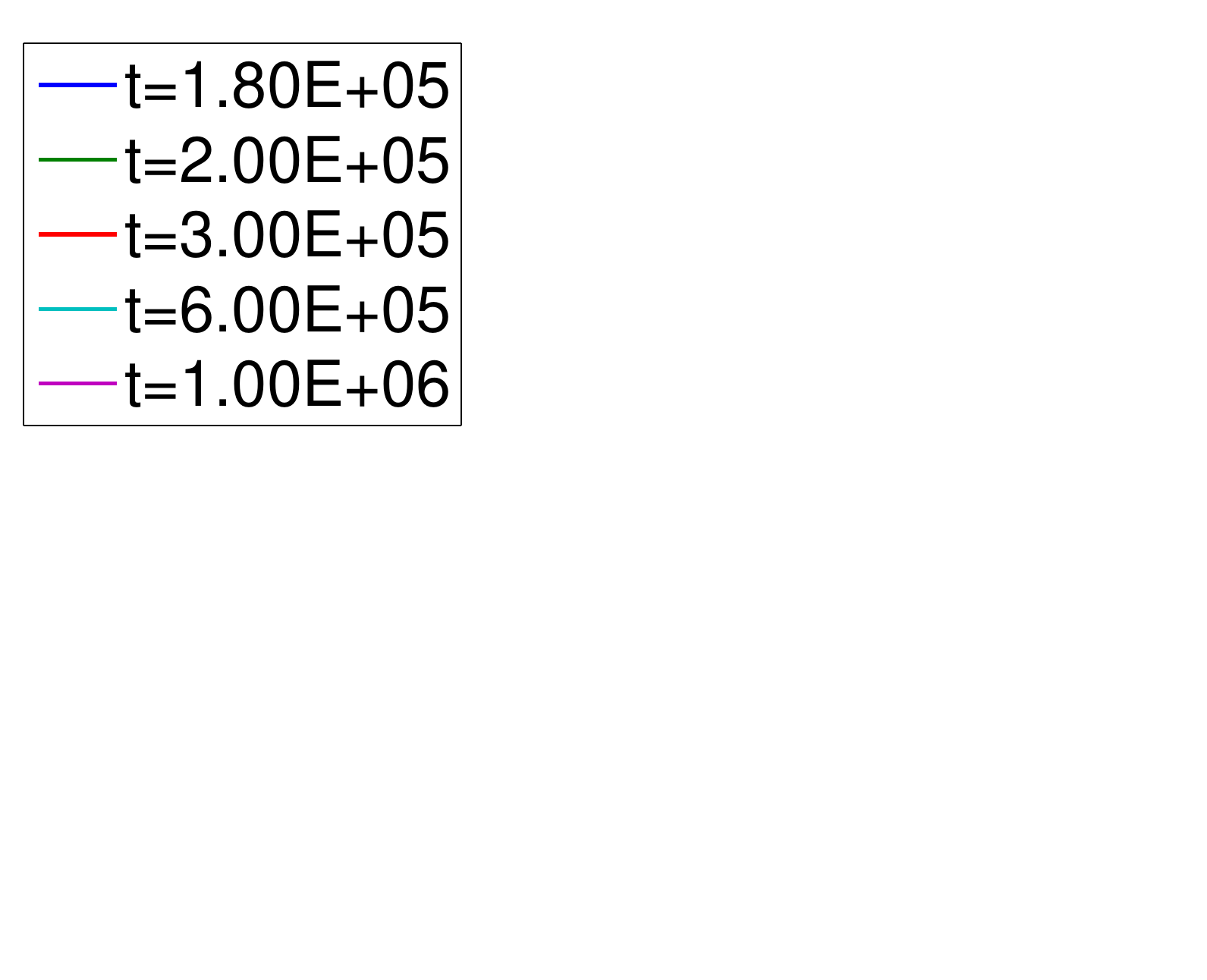}
    \caption{Numerical test case number 4, $L_x=1$ m, $L_1=0.5$ m: Time evolution of the the dissolved hydrogen
    molar density ($\rho_l^h/M^h$) (top right), $p_l$ (top left)
    and $S_g$ (bottom) profiles; during the five last time steps.}
    \label{fig:ct4_2}
\end{figure}

The space discretization step was taken constant equal to $2\cdot
10^{-3}$ m and the time step was  variable, going from $0.33~s$  in
the beginning of the simulation to $16.7\cdot 10^3 \;s $ at the end
of the simulation (see Table~ \ref{tab:NUM}).
Figures~\ref{fig:ct4_1} and~\ref{fig:ct4_2} represent the liquid
pressure $p_l$, the dissolved hydrogen molar density ( equal to
$\rho_l^h/M^h$) and the gas saturation $S_g$ profiles at different
times.

There are essentially two steps:
\begin{itemize}
    \item For $0<t<1.92\cdot 10^5\;s$ (see Figure~\ref{fig:ct4_1}), the initial  gas saturation
    jump moves from $x=0.5~m$, at $t=0$ and reaches $\Gamma_{in}$, the left domain
    boundary, at  $t=1.92\cdot 10^5~s$. During this movement, the
    saturation jump height (initially $\approx0.16\;$) decreases, until
    approximately $0.03$, when it reaches the left
    boundary $\Gamma_{in}$. In front of this discontinuity there is a liquid saturated
        zone, $S_g=0$, and in this zone both the liquid pressure and the hydrogen density are spatially uniform
        (see Figure~\ref{fig:ct4_1}, top ).
        But, while the hydrogen density  remains constant and equal to its initial
        value,
      the liquid pressure  becomes immediately continuous and starts
      growing quickly
       (for instance, $p_l(t=10^3 \; s) \approx1.6\cdot 10^6 \; $Pa), and then more slowly until $t=1.3\cdot 10^5 \;s
       $, when it starts  to slightly decrease.

        In Figure~\ref{fig:ct4_1},  located on the gas saturation
       discontinuity, there are both a high contrast in the  dissolved hydrogen
       concentration (this concentration stays however
       continuous, but with a strong gradient, as seen in the top right of Figure~\ref{fig:ct4_1}), and a
        discontinuity in the liquid pressure
       gradient (see the top left of Figure~\ref{fig:ct4_1}).

      \item For $1.92 \cdot 10^5\;s<t<10^6\;s=T_{fin}$ (see Figure ~\ref{fig:ct4_2}), all the entire domain is now
       unsaturated ($S_g > $0).
          The liquid pressure, the hydrogen density  and the gas saturation profiles are all strictly
          monotonous and continuous, going towards a  spatially uniform distribution,
          corresponding to
          the stationary state(see Figure~\ref{fig:ct4_2}).
 \end{itemize}

 As expected, the system initially out of equilibrium (discontinuity of the gas
 pressure), becomes immediately again in equilibrium (the gas
 pressure is continuous)and
 evolves towards a uniform stationary state (due to the no mass inflow and outflow boundary conditions).
 Although the liquid pressure and the  dissolved hydrogen density are immediately  again continuous for $t > 0 \;
 $, the  hydrogen density still have a locally very strong gradient until $t=1.92\cdot 10^5 \;s $.

At first, and at the very begining($\approx~10^2 \;s $), see top
left of Figure~\ref{fig:ct4_1}, only the liquid pressure evolves in
the liquid saturated zone. Due to a  gas pressure in the unsaturated
zone higher than in the liquid saturated zone ($S_g=0$;
$p_{g}=2.5MPa>p_{l}=1MPa$, for the initial state in
Table\ref{tab:ct4}), and due to the no flow condition imposed on
$\Gamma_{in}$, the liquid in the saturated zone is compressed by the
gas from the unsaturated zone. Then, a liquid gradient pressure
appears around the saturation front and makes the liquid to flow
from the liquid saturated zone towards the unsaturated one, and then
the gas saturation front to move in the opposite direction.

 The very strong hydrogen density gradient (until
$t=1.92\cdot 10^5 \;s
 $), located on the saturation front,
  is due to the competition between the diffusion  and the convective flux
of the dissolved hydrogen around the saturation front: the water
flow convecting the dissolved hydrogen, from left to right, cancels
the smoothing effect of the gas diffusion propagation in the
opposite direction. On the one hand the diffusion is supposed to
reduce the hydrogen concentration contrast, by creating a flux going
from strong concentrations (in the unsaturated zone) towards the low
concentrations (in the liquid saturated zone), and on the other hand
the flow of the liquid phase goes in the opposite direction (left to
right, from $S_g=0$ to $S_g>0$).
 Once the  disequilibrium has
disappeared, the system tends to reach a  uniform stationary state
determined by the mass conservation
 of each component present in the initial state (the system is isolated, with no flow on any of the
boundaries).

\section{Concluding remarks}
\label{Sec: conclusion} From balance equations, constitutive
relations and equations of state, assuming thermodynamical
equilibrium, we have derived a model for describing underground gas
migration in water saturated or unsaturated porous media, including
diffusion of components in phases and capillary effects. In the
second part, we have presented a group of numerical test cases
synthesizing the main challenges concerning  gas migration in a deep
geological repository. These numerical simulations, are based on
simplified but typical situations  in underground nuclear waste
management; they show evidence of the model ability to describe the
gas (hydrogen)  migration, and to treat the difficult problem of
correctly following the saturated and unsaturated regions created by
the gas generation.

\begin{acknowledgements}
This work was partially supported by the GNR MoMaS\\ (PACEN/CNRS,
ANDRA, BRGM, CEA, EDF, IRSN). Most of the work on this paper was
done when Mladen Jurak was visiting, at Universit\'{e} Lyon 1, the
 CNRS-UMR 5208 ICJ.
\end{acknowledgements}

\end{document}